\documentclass{article}

\usepackage[T1]{fontenc}
\usepackage[utf8]{inputenc}
\usepackage[a4paper, left=1in, right=1in, top=1in, bottom=1in]{geometry}  
\usepackage{stmaryrd}
\usepackage{amssymb}
\usepackage{amsmath} 
\usepackage{fancyhdr}

\usepackage[
	backend=bibtex,
	style=numeric,
]{biblatex}
\usepackage[colorlinks=true,citecolor=blue,linkcolor=black]{hyperref}

\usepackage[pdftex]{graphicx} 
\usepackage{algorithm}
\usepackage[noend]{algpseudocode}
\usepackage{pgf}
\usepackage{caption}
\usepackage{subcaption}
\usepackage{bbm}
\usepackage{framed}
\usepackage[most]{tcolorbox}
\usepackage{xcolor}
\definecolor{darkpurple}{RGB}{75, 0, 130} 
\definecolor{darkgreen}{RGB}{0, 100, 0} 
\definecolor{darkorange}{RGB}{255, 140, 0} 

\usepackage{multirow}
\usepackage{longtable}
\usepackage{supertabular}
\usepackage{multicol}

\newcommand{\revision}[1]{{\color{black}#1}} 
\newcommand{\amal}[1]{{\color{black}#1}}

\newcommand{\imad}[1]{{\color{black}#1}}

\title{Computational investigations of a multi-class traffic flow model: mean-field and microscopic dynamics}
\author{A. Habbal \and I. Kissami \and A. Machtalay \and A. Ratnani}
\author{A. Machtalay$^{2,1}$, A. Habbal$^{1,2,*}$, A. Ratnani$^{2}$, I. Kissami$^{3}$
	\\
	\small{\textit{$^1$Universit{\'e} C{\^o}te d'Azur, Inria, CNRS, LJAD UMR 7351, Parc Valrose 06108 NICE CEDEX 2, France} 
	} 
	\\
 \small{\textit{$^2$The UM6P Vanguard Center, Mohammed VI Polytechnic University, Lot 660, Ben Guerir 43150, Morocco} 
	}
 \\
 \small{\textit{$^3$College of Computing, Mohammed VI Polytechnic University, Lot 660, Ben Guerir 43150, Morocco} 
	}
	\\
 	\small{\textit{$^*$Corresponding author: Abderrahmane.Habbal@univ-cotedazur.fr} 
	} 
}
\date{}  

\begin{document}
\pagestyle{fancy}
\fancyhf{}
\fancyhead[L]{\textcolor{magenta}{For citation, please use the following reference: \href{https://doi.org/10.1016/j.trb.2025.103196}{https://doi.org/10.1016/j.trb.2025.103196}}}

\maketitle

\thispagestyle{fancy} 
\pagestyle{plain} 

\begin{abstract}
    We address a multi-class traffic model, for which we computationally assess the ability of mean-field games (MFGs) to yield approximate Nash equilibria for traffic flow games of intractable large finite-players. We introduce \textit{ad hoc} numerical methodologies, with recourse to techniques such as \revision{High-Performance Computing (HPC)} and regularization of \revision{Loose Generalized Minimal Residual (LGMRES)} solvers. The developed apparatus allows us to perform simulations at significantly larger space and time discretization scales. For three generic scenarios of cars and trucks, and three cost functionals, we provide numerous numerical results related to the autonomous vehicles (AVs) traffic dynamics, which corroborate for the multi-class case the effectiveness of the approach emphasized in \cite{huang_game-theoretic_2020}. We additionally provide several original comparisons of macroscopic Nash mean-field speeds with their microscopic versions, allowing us to computationally validate the so-called $\epsilon-$Nash approximation, with a rate slightly better than theoretically expected.
\end{abstract}

\section{Introduction}\label{Intro}


The increasing world movement toward autonomous vehicles (AVs for short) has attracted tremendous research and development effort from researchers and engineers, which opens the way for new significant studies in this field \cite{rahman2023impacts}. 

Paper \cite{huang_game-theoretic_2020} makes a significant contribution to the field of AVs by introducing a game theoretical framework that bridges the gap between microscopic differential games and macroscopic mean field games (MFG for short). This latter, introduced by Lasry and Lions in \cite{lasry_mean_2007} and Huang et al. in \cite{huang_large_2006}, aims to study complex systems in which many agents interact strategically in an evolving stochastic environment. In \cite{huang_game-theoretic_2020}, the authors focus on a single-class of AVs, where all vehicles are autonomous and have similar characteristics. 

However, real life traffic involves many different classes of vehicles, such as human-driven motorcycles, cars and trucks, as well as autonomous vehicles. \revision{To account for such situations, sound models tailored to handle specificities of multiple classes, notably of Autonomous Vehicles are needed.}\\

\revision{
In multi-modal transportation management, the suitable theory for capturing how individual vehicles within a class optimize their strategies while accounting for interactions both within their class and with vehicles from other classes is that of multi-class (a.k.a multi-population) mean-field games  \cite{cirant2015multi}. In some situations, strategic decisions may be of cooperative nature either at the individual scale in each class, or between different classes, leading to optimal control-type formulations.
In our present case, we formulate the multi-class traffic at the macroscopic scale as a Nash Mean-field Game (NMFG for short) following the multi-population classifications in \cite{bensoussan_mean_2018}. NMFG concept is aimed at modeling macroscopic limit behavior of situations where all the vehicles of the same class and in different classes as well interact non-cooperatively. 
}\\

The challenges regarding numerical solution of \revision{both single and multi-class} MFG and NMFG models are two-fold. Firstly, these forward-backward PDE systems involve non-standard and harsh coupling between the unknowns. Secondly, the uniqueness of the numerical solutions is not guaranteed since the Lasry-lions monotonicity condition does not always hold. Several numerical methods are proposed to solve MFG problems. See \cite{lauriere_numerical_2021} for an overview. As in our work, we use a \revision{space-time} method based on finite difference \cite{Achdou2020MeanFG2} and Newton's approach \cite{Achdou2012MeanFG}. The main advantage of this method is that\revision{, for its implementation, } there are no requirements for the separability of the cost function nor the time horizon. However, it may fail to converge in deterministic cases or when the solution is non-smooth. Moreover, it requires additional expensive computations.\\

Our goal is to lead large-scale numerical investigations on an extended multi-class traffic model of one class approach proposed in \cite{huang_game-theoretic_2020}. First, we set up the multi-class traffic framework. Then, we develop original computational techniques to design an efficient and optimized algorithm to solve both one-class MFG and two-class NMFG problems. We provide extensive numerical results of the problem both from large computational scales and from the diversity of scenarios viewpoints.

\subsection*{Paper contribution}

Starting from the theoretical developments regarding the one-class micro-macro and macro-micro bridges of \cite{huang_game-theoretic_2020}, we extend and provide original contributions in several respects. \revision{From the modeling part, the original single-class single-lane framework is extended to a multi-class single-lane traffic model, adapting to multi-class traffic requirements and assumptions (leaving it for future work to address multi-lane and collision avoidance extensions).} 

We have as well investigated three different scenarii regarding initial two-class, cars and trucks, layout. From the computational methodological viewpoint, we use LGMRES with viscous regularization, in a HPC environment, allowing us to address large scale simulations, namely solving several loops of non-linear systems of  up to 88,000,000 unknowns. 

We provide several original comparisons of macroscopic Nash mean-field speeds with their microscopic versions, and for cost functionals as well, allowing us to computationally validate the so-called $\epsilon-$Nash approximation, with  rates slightly better than theoretically expected.

\begin{small}
\begin{center}
\begin{longtable}[!htb]{@{} |ll @{} |ll|l @{}}
\caption{\centering Nomenclature.} 
\label{tab_nomenclature}
\\ \cline{1-4}
$j$            & The $j$-th class of vehicles.         &  $i_j$         & The $i$-th vehicle in the $j$-th class. &  \\ 

$I_j$          & Set of vehicle indexes in class $j$.           &  $N_j$          & Number of vehicles in class $j$.               &  \\ 

$N$     & Total number of vehicles ($\sum_{j}^{}N_j$).    &
$l_j$          & Average length of vehicles in class $j$.
             &  \\ 

$L_j$          & Road section of the $j$-th class.      & 
$\alpha_j$     & Amount of road sections for class $j$.    &  \\ 

$L$     & Total road length ($\sum_{j}^{}L_j \alpha_j$).     &  $u^j_{max}$      & Free flow speed of the $j$-th class.                            &  \\ 

$\rho^j_{jam}$   & Jam density of the $j$-th class.                                &  $v_{i_j}$            & Micro velocity of the $i_j$-th vehicle.          &  \\ 

$x_{i_j}$            & 
Position of the $i_j$-th vehicle.          &  $T$            & Time horizon.         &  \\ 

$u^j$            & Macro velocity field of class $j$.                       &  $\rho^j$         & Density distribution of class $j$.                            &  \\ 

$K$         & Kernel density function.                            &  $\sigma$         & Kernel function bandwidth.                            &  \\ 

$V^j$ & Optimal cost for class $j$.          &  $p_j$ & Gradient of $V^j$ ($\partial_x V^j$).           &  \\ 

$(\rho^*,u^*)$ & Mean-field equilibrium.                      &  $\tilde{v}$    & $\varepsilon$-Nash Equilibrium.              &  \\ 

$\hat{v}$      & Mean field constructed control.              &  $\bar{v}$      & Best response strategy.                      &  \\ 

$f$            & Running cost.                                &  $G$            & Terminal cost.                               &  \\ 

$J$            & Cost functional.                             &  
$H_j$          & Hamiltonian in the $j$-th class.                    &  \\ 

$H^*_j$          & Derivative of $H_j$ w.r.t $p_j$.                    &  w.r.t         & With respect to.                               &  \\ 

$N_t$         & Number of points in time.                               &  $N_x$         & Number of points in space.                               &  \\ 

$\Delta t$         & Time step.                               &  $\Delta x$         & Space step.                               &  \\ 

$(n,k)$         & Time-space index.                               &  $\nu$         & Regularization coefficient.                               &  \\ 

$\hat{\varepsilon}$  & Accuracy.                                    &  MeanRA         & Mean accuracy.                               &  \\ 

MaxRA          & Maximum accuracy.                            &  LWR          & Lighthill-Whitham-Richards model.                            &  \\ 

GLWR          & Generalized LWR.                            &  GS         & Generalized Separable.                            &  \\ 

GNS          & Generalized Non-Separable.                            &  TC          &  Fully segregated configuration with            &  \\ 

CT          & Fully segregated configuration     &  &  cars in front.   &  \\ 

   &   with trucks in front.               &  TCT          & Interlaced configuration with                             &  \\ 

NMFE         & Nash Mean Field Equilibrium.                                      &  & alternating trucks and cars. &   \\ \cline{1-4}

\end{longtable}
\end{center}
\end{small}

\noindent
The organization of the paper is as follows: In section \ref{section2}, we model the multi-class traffic problem at both macro and micro levels and establish the connection between them. In section \ref{section3}, we explore crucial computational methods for designing our algorithm. Finally, in section \ref{section4}, we showcase numerical results.

\section{Multi-class traffic model}\label{section2}
We consider $J$ different classes of vehicles on a single-lane road. \\

\revision{Our microscopic dynamics, see system (\ref{eq.micro.dynamic.system}), do not indeed account for collision avoidance. Incorporating collision avoidance at the microscopic scale in modeling involves considering
 the mechanisms by which vehicles react to potential collisions. To derive an ad hoc multi-class MFG, one may use microscopic models such as the popular Car-Following ones, where inter-vehicle spacing is introduced to prevent rear-end collisions, with dynamics based on reaction times and safe distance rules. One step further would be to consider multi-lane frameworks as in  \cite{festa2018mean}, and incorporate in the lane-changing algorithms constraints or costs for lane changes to avoid conflicts with nearby vehicles. While we presently keep focus on a simple single-lane multi-class model, taking into account multi-lane and collision avoidance extensions is of utmost importance in traffic management, and are sought for future developments.}\\
 
From now on, we denote by $j\in \left \{1,2,...,J\right \}$ the $j$-th class and $i_{j}$ the $i$-th vehicle of the $j$-th class.\\
Every class $j$ \revision{has} $N_{j}$ vehicle, each indexed by $i_{j}\in \mathcal{I}_j=\left \{1+ \sum_{k=1}^{j-1}N_k,...,\sum_{k=1}^{j}N_k \right \}$, where $N=\sum_{j=1}^{J} N_{j}$ the total number of vehicles.\\
\\
\amal{Let denote $u_{max}^{j}$, $\rho^{j}_{jam}$ and $l_{j}$ respectively the free flow speed, the jam density and the average length of the vehicles of the $j$-th class. In a single-lane multi-class traffic scenario, each class of vehicles takes up a different section of the road that can be divided into multiple portions (due to the alternating distribution of vehicle classes) within the total road length $L$. To simplify the representation, we assume that each class $j$ is spread evenly across portions of equal length, referred to as one road section $L_j$. The number of such portions is denoted by $\alpha_j$. Then, the total road length $L=\sum_{j=1}^{J}\alpha_j L_j$.}\\
\\
An additional relation needs to be fulfilled for $j\in \left \{1,2,...,J\right \}$:
\begin{equation}\label{eq.density.lenth.relation}
    \rho^{j}_{jam} =\frac{L_{j}}{l_{j}}
\end{equation}
The expression \eqref{eq.density.lenth.relation} ensures that the total length of vehicles that can fit into the $j$-th road section is exactly the length of the $j$-th road section.\\ 
\\
In the following, we have recourse to a common notation used in game theory for any vector $w$ (with no effective permutation): 
\begin{equation}\label{eq.game.notation}
    w=(w^{j},w^{-j})=((w_{i_{j}},w_{-i_{j}}),(w_{i_{-j}},w_{-i_{-j}}))
\end{equation}
where $j$ denotes the $j$-th class of vehicles, which is a combination of the $i$-th vehicle denoted by $i_j$, and all the other vehicles denoted by $-i_j$. Thus, we can write $w^{j}=(w_{i_{j}},w_{-i_{j}})$.\\ 
As well as $-j$ denotes all classes except the $j$-th class, which is a combination of the $i$-th vehicle denoted by $i_{-j}$, and all the other vehicles denoted by $-i_{-j}$. Thus, we can write $w^{-j}=(w_{i_{-j}},w_{-i_{-j}})$.

\subsection{From microscopic dynamics to the Nash mean-field game model}

In this section, we recall the mean field game approach introduced in \cite{huang_game-theoretic_2020} for one-class traffic model to build the connection between the microscopic and the macroscopic scales, which are based, respectively, on the Lagrangian and the Eulerian observations. We extend their approach to the multi-class traffic model. We assume in all what follows that $j\in \left \{1,2,...,J\right \}$ and $i_{j}\in \mathcal{I}_j$.\\

\noindent {\bf Lagrangian approach.} At the microscopic scale, we follow each vehicle ($i$ of class $j$) position $x_{i_{j}}(t)$ and velocity $v_{i_{j}}(t)$ at time $t \in [0,T]$, where $T>0$ is a fixed time horizon. Thus, the $i_{j}$-th vehicle's dynamic system is the following:
\begin{equation}\label{eq.micro.dynamic.system}
    \left\{\begin{matrix}
\dot{x}_{i_{j}}=v_{i_{j}}(t)\\
x_{i_{j}}(0)=x_{i_{j},0}
\end{matrix}\right.
\end{equation}
where $x_{i_{j},0} \in [0,L]$ is a given initial position.\\
\\
%
Every vehicle $i_j$ tries to minimize its own individual cost functional:
\begin{equation}\label{eq.micro.cost.function}
    J_{i_{j}}(v_{i_{j}},v_{-i_{j}},v^{-j})=\int_{0}^{T}f_{i_{j}}\left (v_{i_{j}}(t),x_{i_{j}}(t),x_{-i_{j}}(t),x^{-j}(t)  \right )dt+G_{i_{j}}(x_{i_{j}}(T))
\end{equation}
where $f_{i_{j}}$ is the running cost function and $G_{i_{j}}$ the terminal cost, accounting for the $i_{j}$-th vehicle's objectives respectively in the current and in the final positions.\\
\\
Our model does not take into account possible collisions at the microscopic scale. Other running cost functions could be considered which account for collisions. One way to do so could be to introduce a velocity-dependent safety distance, see \cite{kamal2014smart}.\\
\\
\revision{We seek for the velocity strategies} $v^{*}=(v^{1,*},...,v^{J,*})$ as a Nash Equilibrium of the following $N$-vehicle differential game:
\begin{equation}\label{eq.micro.nash.equilib}
    \forall v_{i_{j}} \in [0,u^{j}_{max}], \: J_{i_{j}}(v^{*}_{i_{j}},v^{*}_{-i_{j}},v^{-j,*}) \leq J_{i_{j}}(v_{i_{j}},v^{*}_{-i_{j}},v^{-j,*})
\end{equation}
\\
\noindent {\bf Eulerian approach.} At the macroscopic scale, we look at the traffic as a multi-phase flow, and we describe its evolution using an Eulerian observation.\\
We consider the velocity fields $u^{j}(t,x)$ described by:
\begin{equation}\label{eq.macro.velocity}
    u^{j}\left ( t, x_{i_{j}}(t) \right )=v_{i_{j}}(t)
\end{equation}
And $\rho^{j}(t,x)$ the density functions, approximated using a kernel density estimation (see \cite{wkeglarczyk2018kernel}) as follows: 
\begin{equation}\label{eq.macro.KDE}
    \rho^{j}(t,x) \approx \frac{1}{\sigma ~ N_{j}}\sum_{i_{j}\in \mathcal{I}_j}^{} K\left ( \frac{x-x_{i_{j}}(t)}{\sigma} \right )
\end{equation}
where $K$ is a non-negative smoothing kernel function, and $\sigma >0$ is an appropriate bandwidth. In the numerical simulations, we take $K(s)=\frac{1}{\sqrt{2\pi}} e^{-\frac{s^2}{2}}$ and $\sigma=0.05 N_j$.
\\
\\
Since all the vehicles are similar in a same given class $j$, they have the same profile of cost functions. Thus, for all $i_{j}\in \mathcal{I}_j$, we set $f_{j}=f_{i_{j}}$ and $G_{j}=G_{i_{j}}$. We consider a generic vehicle indexed by $j$ (i.e. a typical vehicle in class $j$), then \eqref{eq.micro.cost.function} reads:
\begin{equation}\label{eq.macro.cost.function}
    J_{j}(u^j)=\int_{0}^{T}f_{j}\left (u^j(t,x(t)),\rho(t,x(t))  \right )dt+G_{j}(x(T))
\end{equation}
with 
\begin{equation}\label{eq.generic.dynamic.system}
    \left\{\begin{matrix}
\dot{x}=u^j(t,x(t))\\
x(0)=x_{j,0}
\end{matrix}\right.
\end{equation}
where $x_{j,0} \in [0,L]$ is a given initial position of the generic vehicle $j$, and $\rho=(\rho^1,...\rho^J)$.
\\
Specific types of cost functions $f_j$ and $G_j$ will be addressed in the following sections.\\
\\
\noindent {\bf NMFG System.} Since we see multi-class traffic flow as a multi-phase fluid flow, the evolution of the traffic flow can be described by the following continuity equations, with $j\in \left \{1,2,...,J \right \}$:
\begin{equation}\label{eq.macro.CE}
    \partial _{t}\rho^{j}(t,x)+\partial _{x}(\rho^{j}(t,x) u^{j}(t,x))=0, \quad t\in \left [0,T\right ], x\in \left [0,L\right ]
\end{equation}
Hence, we seek an equilibrium $(\rho^{*},u^{*})$ satisfying the two following conditions:
\begin{enumerate}
    \item $u^{*}=(u^{1,*},...,u^{J,*})$ is a Nash Equilibrium of \eqref{eq.macro.cost.function}-\eqref{eq.generic.dynamic.system}.
    \item $\rho^{*}=(\rho^{1,*},...,\rho^{J,*})$ is a solution to \eqref{eq.macro.CE}, controlled by $u^{*}$.
\end{enumerate}
We introduce the optimal cost, which plays the role of the Bellman value function: \revision{$$V^j(t,x)=\underset{0 \leq u^{j} \leq u^{j}_{max}}{min} \int_{t}^{T}f_{j}\left (u^j(s,x(s)),\rho(s,x(s))  \right )ds+G_{j}(x(T))$$}
By applying the dynamic programming principle of Bellman to  $V^j(t,x)$, we end up with \amal{the equations (ii) and (iii) in} the following  PDE system:
\begin{equation}\label{eq.macro.NMFG.system}
[\text{NMFG}]:    \left\{\begin{matrix}
\forall (t,x) \in \left [0,T\right ] \times \left [0,L\right ]\\
\partial _{t}\rho^{j}(t,x)+\partial _{x}(\rho^{j}(t,x)u^{j}(t,x))=0, \quad (i)\\
\partial_{t}V^{j}(t,x)+H_{j}\left (\partial_{x}V^{j}(t,x), \rho(t,x) \right ) =0, \quad (ii)\\
u^{j}(t,x)= H^{*}_{j}\left (\partial_{x}V^{j}(t,x), \rho(t,x) \right ), \quad (iii)
\end{matrix}\right.
\end{equation}
where $H_{j}$ is the Hamiltonian function for the class $j$ defined by: 
\begin{equation}\label{eq.macro.hamiltonian}
    H_{j}(p_{j}, \rho)=\underset{0 \leq \alpha \leq u^{j}_{max}}{min} \left \{ f_{j}(\alpha, \rho)+\alpha p_{j} \right \} 
\end{equation}
$p_{j}=\partial_{x} V^{j}$ is the derivative of $V^j$ w.r.t $x$, and $H^{*}_{j}=\partial_{p_{j}}H_{j}$ denotes the derivative of $H_{j}$ w.r.t $p_{j}$ yields the optimal velocity field $u^{j}$.\\
Our calculations to obtain \amal{(ii) and (iii) in} \eqref{eq.macro.NMFG.system} follow the classical steps in deriving \revision{Hamilton-Jacobi-Bellman (HJB)} equations (see e.g. \cite{bensoussan_mean_2018})
\\
\\
(i) is the \revision{continuity equation (CE)} known in MFG as the \revision{Fokker-Planck-Kolmogrov (FPK)} equation; this equation is forward in time with an initial condition. (ii) is the backward \revision{HJB} obtained via the dynamic programming principle of Bellman \cite{Bensoussan1984}, with a terminal condition. This forward-backward system is non-linear and coupled via the \textit{feedback-law} (iii). \\
\\
We consider the following initial and terminal conditions:
\begin{equation}\label{eq.macro.BC}
    \left\{\begin{matrix}
\rho^{j}(0,x)=\zeta^j_0(x) & x \in [0,L]\\ G_{j}(x(T))=V_{T}(x) & x \in [0,L]\\
\end{matrix}\right.
\end{equation}
where $V_{T}$ describes the vehicle's preferences on the final position. In the sequel, we assume that there are no preferences. Then, for all $x \in [0,L]$, $V_{T}(x)= 0$.\\
\\
We consider a closed ring road, and use periodic boundary conditions for $\rho$ and $V$ as follow:
\begin{equation}\label{eq.periodic.BC}
\left\{\begin{matrix}
\rho^j(t,0)=\rho^j(t,L)\\ V^j(t,0)=V^j(t,L)\\
\end{matrix}\right.
\end{equation}\\
\amal{In realistic scenarios, speeds are naturally limited. As vehicles cannot move at negative speeds, we have a lower limit of $0$, which means that in the extreme case of traffic jams, vehicles will stop. To ensure that vehicles in each class do not exceed the maximum permitted speed (the free-flow speed that the vehicle can reach when the road is empty), we assume an upper limit of $u_j^{max}$ as a safety constraint. This limit needs to be modeled to reflect traffic dynamics, since in some cases (separable and non-separable), vehicles can choose arbitrarily high or negative speeds to minimize their cost function.}
To take into account the bound constraint on the controls, we use the following projection on $\left [0,u_{max}^j\right ]$ for the (iii) in \eqref{eq.macro.NMFG.system}
\begin{equation}\label{eq.macro.u_condition}
    max\left\{min\left\{u^j(t,x),u^j_{max} \right\},0 \right\}
\end{equation}

\subsection{From the Nash mean-field game equilibrium back to the microscopic velocity controls}\label{sec2.2}

In this section, we extend the approach introduced in \cite{huang_game-theoretic_2020} from one-class to multi-class traffic case, in order to approximate the vehicle's discrete velocity controls from the macroscopic scale back to the microscopic scale, and then we compare it with the classical one in which we directly compute the controls as an $\varepsilon$-Nash equilibrium of the microscopic differential game \cite{Carmona2012ProbabilisticAO}.\\

\noindent {\bf NMFG-Constructed controls.} Our aim is to construct the microscopic discrete optimal controls $(\hat{v}_{i_{j}})_{i_{j}\in \mathcal{I}_j, j\in \left \{1,...,J\right \}}$ from the macroscopic NMFG solutions $(u^{1,*},...,u^{J,*})$ using \eqref{eq.micro.dynamic.system} and \eqref{eq.macro.velocity}. Hence, solving the following system for $i_{j}\in \mathcal{I}_j$ and $j\in \left \{1,2,...,J\right \}$:
\begin{equation}\label{eq.nmfg.constructed.control}
    \left\{\begin{matrix}
\hat{v}_{i_{j}}(t)=u^{j,*}(t,x_{i_{j}}(t))\\ 
\dot{x}_{i_{j}}= \hat{v}_{i_{j}}(t)\\
x_{i_{j}}(0)=x_{i_{j},0}
\end{matrix}\right.
\end{equation}
where $(x_{i_{j},0})_{i_{j}\in \mathcal{I}_j,j\in \left \{1,...,J\right \}}$ are $N$ random samples generated from the same initial density functions $\zeta^j_0$ used in the macroscopic scale. \\
Now, the question that arises is: are they reasonable solutions for the N-vehicle differential game \amal{in \eqref{eq.micro.nash.equilib}}?\\

\noindent {\bf $\varepsilon$-Nash Equilibrium.} Rather than solving \eqref{eq.micro.nash.equilib} directly, which is difficult when the number $N$ of vehicles is large, we search for the approximations $(\Tilde{v}_{i_{j}}(t))_{i_{j}\in \mathcal{I}_j,j\in \left \{1,...,J\right \}}$ as an $\varepsilon$-Nash Equilibrium for \eqref{eq.micro.nash.equilib}; \\
For $i_{j}\in \mathcal{I}_j$ and $j\in \left \{1,...,J\right \}$:
\begin{equation}\label{eq.epsilon.nash.equilib}
    \forall v_{i_{j}} \in [0,u^{j}_{max}], \: J_{i_{j}}(\Tilde{v}_{i_{j}},\Tilde{v}_{-i_{j}},\Tilde{v}^{-j}) \leq J_{i_{j}}(v_{i_{j}},\Tilde{v}_{-i_{j}},\Tilde{v}^{-j}) +\varepsilon_{i_{j}}
\end{equation}\\
Indeed, we require that the NMFE-constructed controls be the best approximations for \eqref{eq.micro.nash.equilib}, in other words, to be an $\varepsilon$-Nash Equilibrium with a lower bound $\hat{\varepsilon}$, given by:
\begin{equation}\label{eq.epsilon.lower.bound}
\begin{split}
 \hat{\varepsilon}_{i_{j}} &= \underset{v_{i_{j}} \in [0,u^{j}_{max}]}{max}\left\{ J_{i_{j}}(\hat{v}_{i_{j}},\hat{v}_{-i_{j}},\hat{v}^{-j})-J_{i_{j}}(v_{i_{j}},\hat{v}_{-i_{j}},\hat{v}^{-j}) \right\} \\
 &= J_{i_{j}}(\hat{v}_{i_{j}},\hat{v}_{-i_{j}},\hat{v}^{-j})- \underset{v_{i_{j}} \in [0,u^{j}_{max}]}{min}\left\{ J_{i_{j}}(v_{i_{j}},\hat{v}_{-i_{j}},\hat{v}^{-j}) \right\}\\
 &= J_{i_{j}}(\hat{v}_{i_{j}},\hat{v}_{-i_{j}},\hat{v}^{-j})- J_{i_{j}}(\Bar{v}_{i_{j}},\hat{v}_{-i_{j}},\hat{v}^{-j})
\end{split}
\end{equation}\\
\amal{The accuracy $\hat{\varepsilon}_{i_{j}}$ in equation (18) compares the cost $J_{i_{j}}(\hat{v}_{i_{j}},\hat{v}_{-i_{j}},\hat{v}^{-j})$ when the vehicle $i_j$ use its NMFE-costructed control $\hat{v}_{i_{j}}$ given that all the other vehicles are using their NMFE-costructed controls, and the cost $J_{i_{j}}(\Bar{v}_{i_{j}},\hat{v}_{-i_{j}},\hat{v}^{-j})$ when the vehicle $i_j$ use its best response strategy $\Bar{v}_{i_{j}}$ given that all the other vehicles are using their NMFE-costructed controls.
The difference between the two costs tells us how much the vehicle underperforms when using the NMFE-costructed control $\hat{v}_{i_{j}}$ instead of its best response strategy $\Bar{v}_{i_{j}}$.}\\


\amal{The concept of exploitability in the context of MFGs reflects the sub-optimality, meaning how much better a vehicle could perform switching to an optimal strategy. Thus, the accuracy $\hat{\varepsilon}_{i_{j}}$ can be interpreted as a form of exploitability for vehicle $i_j$, measuring how much the vehicle $i_j$ could improve by switching from $\hat{v}_{i_{j}}$ to $\Bar{v}_{i_{j}}$. If $\hat{\varepsilon}_{i_{j}}$ is large, it indicates that the constructed control $\hat{v}_{i_{j}}$ is far from optimal and the vehicle $i_j$ can exploit the system by adjusting to the optimal control.}\\
\\
$\bar{v}_{i_{j}}(t)$ is the $i_{j}$-th vehicle's best response strategy, which means the $i_{j}$-th vehicle's velocity control that minimizes its cost when all the other vehicles move with their NMFE-constructed controls, so that:
\begin{equation}\label{eq.best.response.strategy}
    \bar{v}_{i_{j}}=\underset{v_{i_{j}} \in [0,u^{j}_{max}]}{Argmin}\left\{ J_{i_{j}}(v_{i_{j}},\hat{v}_{-i_{j}},\hat{v}^{-j}) \right\}
\end{equation}
with
\begin{equation}\label{eq.optimal.dynamic}
    \left\{\begin{matrix}
    \dot{x}_{m}= v_{m}(t) & \text{ if } m=i_{j}\\
    \dot{x}_{m}= \hat{v}_{m}(t) &  \text{ if } m \neq i_{j}\\
    x_{m}(0)=x_{m,0}
    \end{matrix}\right.
\end{equation}
\amal{There are several studies on approximate Nash equilibria, particularly in the context of MFGs. However, existing results depend on the specific aspects of the model. Paper \cite{lasry_mean_2007} establishes the theoretical framework but with focus on homogeneous populations. Paper \cite{huang_large_2006} offers a theoretical basis for multi-population systems with many agents but in a stochastic dynamic game. Paper \cite{Carmona2018ProbabilisticTO} covers the probabilistic theory in MFGs with stochastic dynamics. Paper \cite{Bensoussan2014MeanFG} study MFGs with heterogeneity in the model with a "dominating player". Whereas, our model is deterministic and introduces heterogeneity by considering multiple classes of vehicles.}

\section{Numerical Methodologies}\label{section3}
In this section, we first present the numerical schemes used to discretize the NMFG system \eqref{eq.macro.NMFG.system}. Then, we outline the numerical methods used to design an efficient solver for \eqref{eq.macro.NMFG.system}.
\subsection{Finite difference schemes }
We consider $N_{x}+1$ discretized space points $\left\{x_0,..., x_{N_x} \right\}$, and $N_{t}+1$ points in time $\left\{t^0,..., t^{N_t} \right\}$. Given mesh sizes $\Delta t=T/N_{t}$ and $\Delta x=L/N_{x}$ respectively in time and space, we use a uniform grid with $x_{k+1} = x_{k} + \Delta x$ and $t^{n+1} = t^{n} + \Delta t$.\\
\\
The numerical approximations of $\rho^j$, $u^j$ and $V^j$ are denoted respectively 
$$
\left\{\begin{matrix}
\rho_{k}^{j,n}\approx \frac{1}{\Delta x}\int_{x_{k-1}}^{x_{k}}\rho^j(t^{n},s) ds & \text{for } n \in \left\{ 0,...,N_{t} \right\}, \: k \in \left\{1,...,N_{x} \right\} \\
u_{k}^{j,n}\approx \frac{1}{\Delta x}\int_{x_{k-1}}^{x_{k}}u^j(t^{n},s) ds & \text{for } n \in \left\{ 0,...,N_{t}-1 \right\}, \: k \in \left\{1,...,N_{x} \right\} \\
V_{k}^{j,n}\approx V^j(t^{n},x_{k}) & \text{for } n \in \left\{ 0,...,N_{t} \right\}, \: k \in \left\{1,...,N_{x} \right\} \\
\end{matrix}\right.
$$
\\
The first continuity equation (i) of \eqref{eq.macro.NMFG.system} is a non-linear and conservative hyperbolic equation, and to discretize it, we use for every class $j$ the following explicit  Lax-Friedrichs scheme: 
\begin{equation}\label{eq.discrete.CE}
    \frac{\rho^{j,n+1}_{k}-\frac{1}{2}\left ( \rho^{j,n}_{k-1}+\rho^{j,n}_{k+1} \right )}{\Delta t}+ \frac{\rho ^{j,n}_{k+1}u^{j,n}_{k+1}-\rho^{j,n}_{k-1}u^{j,n}_{k-1}}{2 \Delta x}=0
\end{equation}
To guarantee the stability of \eqref{eq.discrete.CE}, we consider the following Courant-Friedrichs-Lewy (CFL) condition:
\amal{
\begin{equation}\label{eq.cfl.condition}
    \Delta t = \frac{c ~ \Delta x}{u_{max}}; \quad u_{max}=\underset{j}{\text{max}}\left\{u_{max}^{j} \right\} 
\end{equation}}\\
where the parameter $c$ is chosen such that $0<c \leq 1$. 
\amal{While the CFL condition is classical when considering classical finite difference schemes, we did not lead theoretical convergence analysis, such as a Von Neumann stability analysis, for the present space-time formulation. We however observed, by computational evidence, that violating the CFL condition yields unstable solutions.} In our numerical simulations, we used $c=0.75$\\
\\
For the second HJB equation (ii) of \eqref{eq.macro.NMFG.system}, we use for each class $j$ the following explicit upwind backward difference scheme:
\begin{equation}\label{eq.discrete.HJB}
    \frac{V^{j,n+1}_{k}-V^{j,n}_{k}}{\Delta t}+H_{j} \left (\frac{V^{j,n+1}_{k+1}-V^{j,n+1}_{k}}{\Delta x}, \rho^{1,n}_{k}, \rho^{2,n}_{k}  \right )=0
\end{equation}
and for regularity reasons (see subsection \ref{reg-cont}), sometimes we need to add artificial viscosity to \eqref{eq.discrete.HJB}, discretized as follow:
\amal{
\begin{equation}\label{eq.discrete.viscosity.term}
    +\frac{\nu}{\Delta x ^{2}}\left ( V^{j,n+1}_{k+1}-2 V^{j,n+1}_{k}+V^{j,n+1}_{k-1} \right )
\end{equation}}
Then, the stability condition for \eqref{eq.discrete.viscosity.term} reads:
\begin{equation}\label{eq.viscos.stability.condition}
    \Delta t = \frac{\beta \Delta x^{2}}{\nu}, \quad \text{where } 0<\beta \leq \frac{1}{2}, \quad \text{and } 0<\nu \leq 1
\end{equation}
We use for every class $j$ the following upwind difference scheme to discretize the feedback-law equation (iii) of \eqref{eq.macro.NMFG.system}:
\begin{equation}\label{eq.discrete.feedback}
    u^{j,n}_{k}-H_{j}^{*} \left (\frac{V^{j,n+1}_{k+1}-V^{j,n+1}_{k}}{\Delta x}, \rho^{1,n}_{k}, \rho^{2,n}_{k}  \right )=0
\end{equation}
The initial and terminal conditions \eqref{eq.macro.BC} are discretized by:
\begin{equation}\label{eq.discrete.BC}
    \left\{\begin{matrix}
\rho^{j,0}_{k}-\frac{1}{\Delta x}\int_{x^{j}_{k-1}}^{x^{j}_{k}}\zeta^{j}_{0}(x) dx =0\\
V^{j,N_{t}}_{k}-V_{T}(x^{j}_{k})=0 
\end{matrix}\right.
\end{equation}
Discrete periodic boundary conditions \eqref{eq.periodic.BC} reads:
\begin{equation}\label{eq.discrete.priodic}
    V^{j,n}_{0}=V^{j,n}_{N_{x}}; \quad \rho^{j,n}_{0}=\rho^{j,n}_{N_{x}}
\end{equation}
\amal{In classical finite-difference schemes, generally implicit schemes have no stability concerns, but are computationally expensive compared to explicit schemes. The latter, on the contrary are simpler and easier to implement, but in turn need to fulfill CFL stability condition. In our case, since we are solving the MFG as all-at-once problem, the computational cost of implicit schemes should not be an issue. We implemented and assessed the two versions, and comparison of the results of both explicit and implicit schemes shows that the number of iterations for convergence is nearly the same.}

\subsection{Newton iterations method }
We reformulate the discrete NMFG system \eqref{eq.macro.NMFG.system} in the following form:
\begin{equation}\label{zero_function}
   F(w)=0 
\end{equation}
where $w$ is a very large vector containing all the unknowns : $$\left\{ \rho_{k}^{j,n}\right\}^{0\leq n \leq N_t}_{1 \leq k \leq N_x}, \: \left\{ u_{k}^{j,n}\right\}^{0\leq n \leq N_t-1}_{1 \leq k \leq N_x}, \: \left\{ V_{k}^{j,n}\right\}^{0\leq n \leq N_t}_{1 \leq k \leq N_x}, \: \text{with } j=1,...,J$$ \\
$F$ encodes equations \eqref{eq.discrete.CE}, \eqref{eq.discrete.HJB}, and \eqref{eq.discrete.feedback}, including the initial and terminal conditions \eqref{eq.discrete.BC}, and taking into account the periodic boundary conditions \eqref{eq.discrete.priodic}. \\
\\
To solve \eqref{zero_function}, we use the Newton–Krylov subspace method named Generalized Minimal Residual method (GMRES) introduced in \cite{knoll_jacobian-free_2004} since it is the most suitable method for a large sparse system with non-linear equations. \\
\\
We first lead computational experiments for the case of one-class traffic proposed in \cite{huang_game-theoretic_2020}, in order to design and validate optimal algorithmic choices. Afterwards, the algorithm is extended to handle the case of two-class traffic.\\

\noindent {\bf LGMRES solver.} "Loose" GMRES solver can be viewed as an acceleration technique for GMRES with or without preconditioning. The LGMRES is one of the GMRES extensions that is more flexible and capable of handling non-normal matrices while maintaining a balance between the convergence rate and computational cost of alternative GMRES methods. Experimental results in \cite{WOS:000230502200005} demonstrate that LGMRES does not require as many iterations as GMRES due to alternating residual vectors.\\
\\
However, in our case, LGMRES does not help much without preconditioning. \revision{We} assess this using the one class [LWR] example proposed in \cite{huang_game-theoretic_2020} (page 14) with the numerical inputs (page 19). As shown in Table \ref{tab_lwr1_approxJ} (see "without preconditioning" column). Indeed, after $1000$ iterations, the residual did not fulfill convergence criteria (Res$\leq6.10^{-6}$).
\begin{table}[H]
    \centering
    \caption{\centering One-class {[}LWR{]}: number of iterations for convergence (iter), final residual (Res), root mean square error (RMSE), and CPU time (in seconds) for different spatial-temporal grid sizes $N_x \times N_t$. The comparison includes the following scenarios: without preconditioning, with preconditioning but without \revision{multi-grid}, and with both preconditioning and \revision{multi-grid} methods. 
    }
        \resizebox{\textwidth}{!}{%
        \begin{tabular}{cc|cc|cc|cccc|l}
        \cline{3-10}
         &  & \multicolumn{2}{c|}{\begin{tabular}[c]{@{}c@{}}without \\ preconditioning\end{tabular}} &  \multicolumn{2}{c|}{\begin{tabular}[c]{@{}c@{}}with \\ preconditioning\end{tabular}} &  \multicolumn{4}{c|}{\begin{tabular}[c]{@{}c@{}}with\\preconditioning \& \revision{multi-grid}\end{tabular}} & \\ \cline{1-10}
        \multicolumn{1}{|c|}{$N_{x}$} & $N_{t}$ & \multicolumn{1}{c|}{iter} & Res & \multicolumn{1}{c|}{iter} & Res & \multicolumn{1}{c|}{iter} & \multicolumn{1}{c|}{Res} & \multicolumn{1}{c|}{RMSE} & time &  \\ \cline{1-10}
        \multicolumn{1}{|c|}{30} & 120 & \multicolumn{1}{c|}{1000} & 4.96 & \multicolumn{1}{c|}{9} & 1.4e-07 & \multicolumn{1}{c|}{2} & \multicolumn{1}{c|}{7.2e-07} & \multicolumn{1}{c|}{0.028} & 3.5 & \\ 
        \multicolumn{1}{|c|}{60} & 240 & \multicolumn{1}{c|}{1000} & 5.42 & \multicolumn{1}{c|}{30} & 2.3e-06 & \multicolumn{1}{c|}{3} & \multicolumn{1}{c|}{1.6e-08} & \multicolumn{1}{c|}{0.022} & 2e01 & \\ 
        \multicolumn{1}{|c|}{120} & 480 & \multicolumn{1}{c|}{1000} & 5.01 & \multicolumn{1}{c|}{117} & 4.5e-06 & \multicolumn{1}{c|}{5} & \multicolumn{1}{c|}{7.5e-07} & \multicolumn{1}{c|}{0.016} & 1e02 & \\ 
        \multicolumn{1}{|c|}{240} & 960 & \multicolumn{1}{c|}{1000} & 5.33 & \multicolumn{1}{c|}{659} & 6.0e-06 & \multicolumn{1}{c|}{15} & \multicolumn{1}{c|}{2.0e-06} & \multicolumn{1}{c|}{0.011} & 2e03 & \\ 
        \multicolumn{1}{|c|}{480} & 1920 & \multicolumn{1}{c|}{1000} & 12.80 & \multicolumn{1}{c|}{1000} & 0.262 & \multicolumn{1}{c|}{168} & \multicolumn{1}{c|}{5.9e-06} & \multicolumn{1}{c|}{0.008} & 7e04 &
        \\ 
        \multicolumn{1}{|c|}{960} & 3840 & \multicolumn{1}{c|}{1000} & 12.01 & \multicolumn{1}{c|}{1000} & 0.404 & \multicolumn{1}{c|}{840} & \multicolumn{1}{c|}{6.0e-06} & \multicolumn{1}{c|}{0.005} &  2e06
        & \\ \cline{1-10}
        \end{tabular}%
    }\label{tab_lwr1_approxJ}
\end{table}

\noindent {\bf Preconditioning.} To accelerate the LGMRES nonlinear solver, we construct a preconditioner to be used at the current iterations by taking the inverse of the Jacobian matrix of $F(w)$ in \eqref{zero_function}. Here, the inverse is computed using the LU decomposition. An approximation of the exact Jacobian is sufficient while dropping the forward-backward coupling parts of $F$, namely the parts with $H_j$ in \eqref{eq.discrete.HJB} and $H_{j}^{*}$ in \eqref{eq.discrete.feedback}. We observe that only one calculation of the preconditioner in the first iteration is needed.\\
\\
The numerical experiment for [LWR] shows in Table \ref{tab_lwr1_approxJ} that the preconditioning strategy favors the solver's convergence. However, the convergence rate slows down for largest $N_x \times N_t$, which is the reason why we have recourse to the \revision{multi-grid} method described in the next paragraph. \\

\noindent {\bf \revision{Multi-grid}.} We remark from the above that using initial guesses set at zero requires more iterations, especially for a fine grid, see "with preconditioning" column in Table \ref{tab_lwr1_approxJ}. The basic idea of the \revision{multi-grid} method is to solve the problem on a fine grid with a precomputed initial guess \cite{WOS:000231357700001}.\\
\\
We start by solving the problem on a coarse grid. Then, using linear cubic interpolation; we interpolate the solution from the coarse to the fine grid, which provides a good initial guess for the fine one. We repeat the procedure as many times as necessary, starting from a grid so that the problem is easy to solve. In the present case, we take for the coarser grid $N_x=15, N_t=60$, with an initial guess set at zero. We observe in the column "with preconditioning and \revision{multi-grid}" in Table \ref{tab_lwr1_approxJ} that the \revision{multi-grid} method leads to a significant reduction in iterations. \\
\\
For both one class traffic cost functions [Separable] and [Non-Separable] proposed in \cite{huang_game-theoretic_2020} (page16) with the numerical inputs (page19), the results in Table \ref{tab_sep1_approxJ.1} and Table \ref{tab_nonsep1_approxJ.1} show that while using an approximate Jacobian for preconditioning, after $1000$ \revision{iterations}, we still are far from convergence, even with preconditioning and \revision{multi-grid} methods. This is likely due to the non-smoothness of the solution. One way to fix this issue is based on regularization-continuation techniques presented in the next paragraph.

\begin{table}[H]
    \centering
    \caption{\centering One-class {[}Separable{]}: regularization coefficient ($\nu$), number of iterations for convergence (iter), final residual (Res), root mean square error (RMSE), and CPU time (in seconds) for different spatial-temporal grid sizes $N_x \times N_t$. Using preconditioning and \revision{multi-grid} methods. Without regularization (a) and with regularization-continuation (b).
    }
    \begin{subtable}{0.30\linewidth}
      \centering
        \caption{\centering $\nu=0.0$.}
        \label{tab_sep1_approxJ.1}
        \resizebox{\textwidth}{!}{%
        \begin{tabular}{|c|c|c|c|l}
        \cline{1-4}
        $N_{x}$ & $N_{t}$ & iter & Res &  \\ \cline{1-4}
        30 & 120 & 1000 & 0.35 &  \\ 
        60 & 240 & 1000 & 0.71 &  \\ 
        120 & 480 & 1000 & 0.9 &  \\ 
        240 & 960 & 1000 & 1.1 &  \\ \cline{1-4}
        \end{tabular}%
        }
    \end{subtable}%
    \hfill 
    \begin{subtable}{0.56\linewidth}
      \centering
        \caption{\centering $0.01 \leq \nu \leq 0.04$.}
        \label{tab_sep1_approxJ.2}
        \resizebox{\textwidth}{!}{%
        \begin{tabular}{|c|c|c|c|c|c|c|l}
        \cline{1-7}
        $N_{x}$ & $N_{t}$ & $\nu$ & iter & Res & RMSE & time & \\ \cline{1-7}
        30 & 240 & 0.04 & 5 & 2.1e-06 & 0.0100 & 2e01 & \\ 
        60 & 720 & 0.03 & 5 & 2.2e-06 & 0.0068 & 1e02 & \\ 
        120 & 1920 & 0.02 & 19 & 3.3e-06 & 0.0107 & 3e03 & \\ 
        240 & 3840 & 0.01 & 959 & 5.5e-06 & 0.0170 & 5e05 & \\ \cline{1-7}
        \end{tabular}%
        }
    \end{subtable}%
     
    \label{tab_sep1_approxJ}
\end{table}


\begin{table}[H]
    \centering
    \caption{\centering One-class {[}Non-Separable{]}: regularization coefficient ($\nu$), number of iterations for convergence (iter), final residual (Res), root mean square error (RMSE), and CPU time (in seconds) for different spatial-temporal grid sizes $N_x \times N_t$. Using preconditioning and \revision{multi-grid} methods. Without regularization (a) and with regularization-continuation (b).
    }
    \begin{subtable}{0.3\linewidth}
      \centering
        \caption{\centering $\nu=0.0$.}
        \label{tab_nonsep1_approxJ.1}
        \resizebox{\textwidth}{!}{%
        \begin{tabular}{|c|c|c|c|l}
        \cline{1-4}
        $N_{x}$ & $N_{t}$ & iter & Res &  \\ \cline{1-4}
        30 & 120 & 1000 & 0.097 &  \\ 
        60 & 240 & 1000 & 0.137 &  \\ 
        120 & 480 & 1000 & 0.164 &  \\ 
        240 & 960 & 1000 & 0.9 &  \\ \cline{1-4}
        \end{tabular}%
        }
    \end{subtable}%
    \hfill
    \begin{subtable}{0.56\linewidth}
      \centering
        \caption{\centering $0.01 \leq \nu \leq 0.04$.}
        \label{tab_nonsep1_approxJ.4}
        \resizebox{\textwidth}{!}{%
        \begin{tabular}{|c|c|c|c|c|c|c|l}
        \cline{1-7}
        $N_{x}$ & $N_{t}$ & $\nu$ & iter & Res & RMSE & time & \\ \cline{1-7}
        30 & 240 & 0.04 & 2 & 4.1e-06 & 0.0083 & 6e01 & \\ 
        60 & 720 & 0.03 & 3 & 1.6e-07 & 0.0057 &  1e02 & \\ 
        120 & 1920 & 0.02 & 5 & 1.0e-06 & 0.0068 & 7e02 & \\ 
        240 & 3840 & 0.01 & 112 & 4.5e-06 & 0.0040 & 5e04 & \\ \cline{1-7}
        \end{tabular}%
        }
    \end{subtable} 
    \label{tab_nonsep1_approxJ}
\end{table}

\noindent {\bf Regularization-continuation.}\label{reg-cont}  An alternative way to fix the non-convergence of solutions for [Separable] and [Non-Separable] while using approximate Jacobian for preconditioning is through regularization techniques as proposed in \cite{achdou2013hamiltonbook}. \\
\\
Indeed, a smoothing effect is introduced by adding the discrete viscosity term \eqref{eq.discrete.viscosity.term} to the discrete HJB equation \eqref{eq.discrete.HJB}, where $ 0\leq \nu \leq 0.05 $ is a regularization parameter. To be much closer to the viscosity-free exact solution of HJB equation, $\nu$ should be as close to zero as possible. Since the Newton method converges better when $\nu$ is large, it is also possible to use the continuation technique proposed in \cite{Achdou2020MeanFG}. We start by solving the problem with a large $\nu$, then use the solution as an initial guess for the fine grid with a smaller $\nu$, and so on until reaching the desired value of $\nu$ close to zero as shown in Table \ref{tab_sep1_approxJ.2} and Table \ref{tab_nonsep1_approxJ.4}. \\
\\
\revision{We employ a continuation technique to approach a desired  $\nu$ value that allows the numerical solution to closely approximate the exact one while mitigating numerical challenges, such as the occurrence of singular Jacobians during Newton’s iterations. Indeed,
when $N_x$ and $N_t$ are large consequently $\Delta x$ and $\Delta t$ are small and when we take $\nu$ close to zero, the viscous term \eqref{eq.discrete.viscosity.term} and the terms involving spatial–temporal variations diminish in the discrete equations, meaning that certain equations in the system begin to lose their independence. Consequently, the Jacobian matrix is no longer full-rank causing a near-singular (or numerically singular) matrix.}\\

\noindent {\bf Root Mean Square Error (RMSE).} 
We have exact data solution $(\omega_{i})_{i=1,...,N}$ which is the solution on a coarser grid and we use it to create, through interpolation, an initial guess $(\tilde{\omega}_{i})_{i=1,...,N}$ on a finer grid.
We use one of the most popular error metrics to assess the difference between exact data solution $(\omega_{i})_{i=1,...,N}$ and predicted one as an interpolated initial guess $(\tilde{\omega}_{i})_{i=1,...,N}$ where $N$ is the number of data points, and $\omega$ is the solution of equation \eqref{zero_function}. 
\begin{equation}\label{eq.rmse}
    \text{RMSE}=\sqrt{\frac{\sum_{i=1}^{N}\left\| \omega_{i}-\tilde{\omega}_{i}\right\|^{2}}{N}}
\end{equation}
If the RMSE decreases with decreasing the grid size, it suggests the efficiency of this approach to accelerate the convergence of solutions and that the method is working correctly, it indicates also that the solution is smooth enough that the interpolation does not introduce important errors.\\

\noindent {\bf Parallel implementation.} 
Another way to fix the non-convergence of solutions for [Separable] and [Non-Separable] consists of using the exact Jacobian matrix for preconditioning to guarantee convergence, as shown in Table \ref{tab_class1_exactJ}. Unfortunately, this does not work for finer grids ($N_x>240$) due to memory issues.
\begin{table}[H]
\centering
\caption{\centering One-class {[}LWR{]}, {[}Separable{]}, and {[}Non-Separable{]}:  number of iterations for convergence (iter), final residual (Res), root mean square error (RMSE), and CPU time (in seconds) for different spatial-temporal grid sizes $N_x \times N_t$. Using preconditioning and \revision{multi-grid} methods but without the regularization-continuation method. 
}
\label{tab_class1_exactJ}
\resizebox{\textwidth}{!}{%
\begin{tabular}{cc|cccc|cccc|cccc|l}
\cline{3-14}
 &  & \multicolumn{4}{c|}{LWR} & \multicolumn{4}{c|}{Separable} & \multicolumn{4}{c|}{Non-Separable} & \multicolumn{1}{c}{} \\ \cline{1-14}
\multicolumn{1}{|c|}{$N_x$} & $N_t$ & \multicolumn{1}{c|}{iter} & \multicolumn{1}{c|}{Res} & \multicolumn{1}{c|}{RMSE} & time & \multicolumn{1}{c|}{iter} & \multicolumn{1}{c|}{Res} & \multicolumn{1}{c|}{RMSE} & time & \multicolumn{1}{c|}{iter} & \multicolumn{1}{c|}{Res} & \multicolumn{1}{c|}{RMSE} & time &  \\ \cline{1-14}
\multicolumn{1}{|c|}{30} & 120 & \multicolumn{1}{c|}{2} & \multicolumn{1}{c|}{5.1e-06} & \multicolumn{1}{c|}{0.028} & 3e01 & \multicolumn{1}{c|}{4} & \multicolumn{1}{c|}{1.9e-08} & \multicolumn{1}{c|}{0.035} & 3e01 & \multicolumn{1}{c|}{4} & \multicolumn{1}{c|}{1.4e-08} & \multicolumn{1}{c|}{0.030} & 2e01 &  \\ 
\multicolumn{1}{|c|}{60} & 240 & \multicolumn{1}{c|}{3} & \multicolumn{1}{c|}{1.8e-09} & \multicolumn{1}{c|}{0.022} & 2e02 & \multicolumn{1}{c|}{5} & \multicolumn{1}{c|}{4.6e-11} & \multicolumn{1}{c|}{0.034} & 2e02 & \multicolumn{1}{c|}{4} & \multicolumn{1}{c|}{1.8e-06} & \multicolumn{1}{c|}{0.025} & 9e01 &  \\ 
\multicolumn{1}{|c|}{120} & 480 & \multicolumn{1}{c|}{3} & \multicolumn{1}{c|}{3.0e-08} & \multicolumn{1}{c|}{0.016} & 7e02 & \multicolumn{1}{c|}{6} & \multicolumn{1}{c|}{2.6e-09} & \multicolumn{1}{c|}{0.027} & 1e03 & \multicolumn{1}{c|}{5} & \multicolumn{1}{c|}{9.1e-11} & \multicolumn{1}{c|}{0.018} & 9e02 &  \\ 
\multicolumn{1}{|c|}{240} & 960 & \multicolumn{1}{c|}{3} & \multicolumn{1}{c|}{1.9e-07} & \multicolumn{1}{c|}{0.011} & 9e03 & \multicolumn{1}{c|}{9} & \multicolumn{1}{c|}{4.8e-06} & \multicolumn{1}{c|}{0.019} & 4e04 & \multicolumn{1}{c|}{5} & \multicolumn{1}{c|}{1.4e-08} & \multicolumn{1}{c|}{0.012} & 5e03 &  \\ 
\multicolumn{1}{|c|}{480} & 1920 & \multicolumn{1}{c|}{3} & \multicolumn{1}{c|}{3.3e-07} & \multicolumn{1}{c|}{0.008} & 7e04 & \multicolumn{1}{c|}{-} & \multicolumn{1}{c|}{-} & \multicolumn{1}{c|}{-} & - & \multicolumn{1}{c|}{-} & \multicolumn{1}{c|}{-} & \multicolumn{1}{c|}{-} & - &  \\ \cline{1-14}
\end{tabular}%
}
\end{table}

In order to improve efficiency and avoid memory issues, we use recent techniques to parallelize and speed up our algorithm. To achieve this, the following three standard Python libraries are used:\\
\\
The first library, mpi4py,  which is updated in \cite{Dalcin2021mpi4pySU}, provides access to the majority of functions of the Message Passing Interface (MPI for short). This latter allows communication between processes running on separate nodes in a parallel system.\\
\\
The second library is the powerful Portable Extensible Toolkit for Scientific Computation (PETSc for short), specifically designed for large-scale numerical simulations in parallel environments and supports different programming languages, including Python, see \cite{abhyankar2018petsc,balay2019petsc}. Further information can be found on the official website \cite{petsc-web-page}.\\
\\
The third library is Pyccel, a Python extension language using accelerators as described in \cite{bourne2023pyccel}. Pyccel is designed to make it easier to write a high-performance algorithm in Python that can be translated to C or Fortran and support the use of accelerators, such as GPUs, to increase performance further.\\
\\
The algorithm to solve the overall NMFG system is implemented in Python. \imad{The numerical experiments were conducted on the TOUBKAL Supercomputer \footnote{https://toubkal.um6p.ma}, featuring a PowerEdge R840 server with four Intel Xeon Platinum 8276 CPUs (112 cores in total) and 1.5 TB of RAM. The CPU governor was set to 'Performance' mode, and jobs were managed using SLURM to optimize resource allocation and prevent memory overuse.}\\
\\
To validate our algorithm numerically, we reproduce in Appendix \ref{appendix_one-class} the results of paper \cite{huang_game-theoretic_2020} for one-class traffic. The performance of the parallel algorithm in terms of CPU time is compared with a sequential one in Table \ref{tab_class1_petsc_lwr}, Table \ref{tab_class1_petsc_separable}, and Table \ref{tab_class1_petsc_nonseparable}. We observe that nothing but running the code in parallel while using only one core reduces the CPU time. We consider $t_{old}$ the CPU time of the sequential code (LGMRES), $t_{new}$ the CPU time of the parallel one \revision{(PETSc)} with only $1$ CPU core, and we compute the percentage improvement of the saved time (PIST) using the following formula:
$$\left (\frac{t_{old}-t_{new}}{t_{old}} \times 100  \right ) \%$$
\begin{table}[H]
\centering
\caption{\centering One-class {[}LWR{]}: number of iterations for convergence (iter), root mean square error (RMSE), and CPU time (in seconds) for different spatial-temporal grid sizes $N_x \times N_t$. Comparison between sequential (LGMRES) and parallel results with different numbers of CPU cores. 
.}
\label{tab_class1_petsc_lwr}
\resizebox{\textwidth}{!}{%
\begin{tabular}{cccc|ccccccccc|c}
\cline{5-13}
                           &                           &                           &      & \multicolumn{9}{c|}{parallel (CPU cores)}                                                                                                                                                                                                                        &                      \\ \cline{3-13}
                           & \multicolumn{1}{c|}{}     & \multicolumn{2}{c|}{LGMRES}       & \multicolumn{2}{c|}{1}                            & \multicolumn{2}{c|}{16}                          & \multicolumn{2}{c|}{56}                          & \multicolumn{2}{c|}{112}                        & \multirow{2}{*}{RMSE} &                      \\ \cline{1-12}
\multicolumn{1}{|c|}{$N_x$}   & \multicolumn{1}{c|}{$N_t$}   & \multicolumn{1}{c|}{iter} & time  & \multicolumn{1}{c|}{iter} & \multicolumn{1}{c|}{time}  & \multicolumn{1}{c|}{iter} & \multicolumn{1}{c|}{time}  & \multicolumn{1}{c|}{iter} & \multicolumn{1}{c|}{time}  & \multicolumn{1}{c|}{iter} & \multicolumn{1}{c|}{time} &                       &                      \\ \cline{1-13}
\multicolumn{1}{|c|}{30}   & \multicolumn{1}{c|}{120}  & \multicolumn{1}{c|}{2}    & 30.4822 & \multicolumn{1}{c|}{4}    & \multicolumn{1}{c|}{0.6479}  & \multicolumn{1}{c|}{4}    & \multicolumn{1}{c|}{0.1957}  & \multicolumn{1}{c|}{4}    & \multicolumn{1}{c|}{0.1865}  & \multicolumn{1}{c|}{4}     & \multicolumn{1}{c|}{0.1614}    & 0.0289                 &                      \\ 
\multicolumn{1}{|c|}{60}   & \multicolumn{1}{c|}{240}  & \multicolumn{1}{c|}{3}    & 173.088 & \multicolumn{1}{c|}{4}    & \multicolumn{1}{c|}{2.9317}  & \multicolumn{1}{c|}{4}    & \multicolumn{1}{c|}{0.8708}  & \multicolumn{1}{c|}{4}    & \multicolumn{1}{c|}{0.7639}  & \multicolumn{1}{c|}{4}     & \multicolumn{1}{c|}{0.6494}    & 0.0226                 &                      \\ 
\multicolumn{1}{|c|}{120}  & \multicolumn{1}{c|}{480}  & \multicolumn{1}{c|}{3}    & 702.323 & \multicolumn{1}{c|}{4}    & \multicolumn{1}{c|}{13.152} & \multicolumn{1}{c|}{4}    & \multicolumn{1}{c|}{3.9964}  & \multicolumn{1}{c|}{4}    & \multicolumn{1}{c|}{3.5443}  & \multicolumn{1}{c|}{4}     & \multicolumn{1}{c|}{2.8693}    & 0.0161                 &                      \\ 
\multicolumn{1}{|c|}{240}  & \multicolumn{1}{c|}{960}  & \multicolumn{1}{c|}{3}    & 9637.26 & \multicolumn{1}{c|}{4}    & \multicolumn{1}{c|}{61.594} & \multicolumn{1}{c|}{4}    & \multicolumn{1}{c|}{17.838} & \multicolumn{1}{c|}{4}    & \multicolumn{1}{c|}{16.601} & \multicolumn{1}{c|}{4}     & \multicolumn{1}{c|}{12.605}    & 0.0114                &                      \\ 
\multicolumn{1}{|c|}{480}  & \multicolumn{1}{c|}{1920} & \multicolumn{1}{c|}{3}    & 72571.4 & \multicolumn{1}{c|}{4}    & \multicolumn{1}{c|}{312.04} & \multicolumn{1}{c|}{4}    & \multicolumn{1}{c|}{81.798} & \multicolumn{1}{c|}{4}    & \multicolumn{1}{c|}{74.985} & \multicolumn{1}{c|}{4}     & \multicolumn{1}{c|}{56.501}    & 0.0080               &                      \\ 
\multicolumn{1}{|c|}{960}  & \multicolumn{1}{c|}{3840} & \multicolumn{1}{c|}{-}    & -    & \multicolumn{1}{c|}{-}    & \multicolumn{1}{c|}{-}    & \multicolumn{1}{c|}{4}    & \multicolumn{1}{c|}{415.64} & \multicolumn{1}{c|}{4}    & \multicolumn{1}{c|}{342.05} & \multicolumn{1}{c|}{4}     & \multicolumn{1}{c|}{265.41}    & 0.0056                 &                      \\ 
\multicolumn{1}{|c|}{1920} & \multicolumn{1}{c|}{7680} & \multicolumn{1}{c|}{-}    & -    & \multicolumn{1}{c|}{-}    & \multicolumn{1}{c|}{-}    & \multicolumn{1}{c|}{-}    & \multicolumn{1}{c|}{-}    & \multicolumn{1}{c|}{4}    & \multicolumn{1}{c|}{2005.9} & \multicolumn{1}{c|}{4}     & \multicolumn{1}{c|}{1374.0}    & 0.0039                 &                      \\ \cline{1-13}
\end{tabular}
}
\end{table}
\begin{table}[H]
\centering
\caption{\centering One-class {[}Separable{]}: number of iterations for convergence (iter), root mean square error (RMSE), and CPU time (in seconds) for different spatial-temporal grid sizes $N_x \times N_t$. Comparison between sequential (LGMRES) and parallel results with different numbers of CPU cores. 
}
\label{tab_class1_petsc_separable}
\resizebox{\textwidth}{!}{%
\begin{tabular}{cccc|ccccccccc|c}
\cline{5-13}
                           &                           &                           &      & \multicolumn{9}{c|}{parallel (CPU cores)}                                                                                                                                                                                                                        &                      \\ \cline{3-13}
                           & \multicolumn{1}{c|}{}     & \multicolumn{2}{c|}{LGMRES}       & \multicolumn{2}{c|}{1}                            & \multicolumn{2}{c|}{16}                          & \multicolumn{2}{c|}{56}                          & \multicolumn{2}{c|}{112}                        & \multirow{2}{*}{RMSE} &                      \\ \cline{1-12}
\multicolumn{1}{|c|}{$N_x$}   & \multicolumn{1}{c|}{$N_t$}   & \multicolumn{1}{c|}{iter} & time  & \multicolumn{1}{c|}{iter} & \multicolumn{1}{c|}{time}  & \multicolumn{1}{c|}{iter} & \multicolumn{1}{c|}{time}  & \multicolumn{1}{c|}{iter} & \multicolumn{1}{c|}{time}  & \multicolumn{1}{c|}{iter} & \multicolumn{1}{c|}{time} &                       &                      \\ \cline{1-13}
\multicolumn{1}{|c|}{30}   & \multicolumn{1}{c|}{120}  & \multicolumn{1}{c|}{4}    & 3.77560 & \multicolumn{1}{c|}{5}    & \multicolumn{1}{c|}{0.8363}  & \multicolumn{1}{c|}{5}    & \multicolumn{1}{c|}{0.2356}  & \multicolumn{1}{c|}{5}    & \multicolumn{1}{c|}{0.21561}  & \multicolumn{1}{c|}{5}     & \multicolumn{1}{c|}{0.23177}    & 0.0204                 &                      \\ 
\multicolumn{1}{|c|}{60}   & \multicolumn{1}{c|}{240}  & \multicolumn{1}{c|}{5}    & 22.8240 & \multicolumn{1}{c|}{5}    & \multicolumn{1}{c|}{4.0176}  & \multicolumn{1}{c|}{5}    & \multicolumn{1}{c|}{1.1025}  & \multicolumn{1}{c|}{5}    & \multicolumn{1}{c|}{0.92241}  & \multicolumn{1}{c|}{5}     & \multicolumn{1}{c|}{0.92102}    & 0.0136                 &                      \\ 
\multicolumn{1}{|c|}{120}  & \multicolumn{1}{c|}{480}  & \multicolumn{1}{c|}{6}    & 134.806 & \multicolumn{1}{c|}{6}    & \multicolumn{1}{c|}{23.654} & \multicolumn{1}{c|}{6}    & \multicolumn{1}{c|}{6.3526}  & \multicolumn{1}{c|}{6}    & \multicolumn{1}{c|}{5.33464}  & \multicolumn{1}{c|}{6}     & \multicolumn{1}{c|}{5.06342}    & 0.0087                 &                      \\ 
\multicolumn{1}{|c|}{240}  & \multicolumn{1}{c|}{960}  & \multicolumn{1}{c|}{9}    & 40000.0 & \multicolumn{1}{c|}{6}    & \multicolumn{1}{c|}{128.77} & \multicolumn{1}{c|}{6}    & \multicolumn{1}{c|}{31.042} & \multicolumn{1}{c|}{6}    & \multicolumn{1}{c|}{26.4452} & \multicolumn{1}{c|}{6}     & \multicolumn{1}{c|}{24.0049}    & 0.0054                 &                      \\ 
\multicolumn{1}{|c|}{480}  & \multicolumn{1}{c|}{1920} & \multicolumn{1}{c|}{-}    & -    & \multicolumn{1}{c|}{7}    & \multicolumn{1}{c|}{924.39} & \multicolumn{1}{c|}{7}    & \multicolumn{1}{c|}{193.27} & \multicolumn{1}{c|}{7}    & \multicolumn{1}{c|}{154.744} & \multicolumn{1}{c|}{7}     & \multicolumn{1}{c|}{114.693}    & 0.0032                 &                      \\ 
\multicolumn{1}{|c|}{960}  & \multicolumn{1}{c|}{3840} & \multicolumn{1}{c|}{-}    & -    & \multicolumn{1}{c|}{-}    & \multicolumn{1}{c|}{-}    & \multicolumn{1}{c|}{9}    & \multicolumn{1}{c|}{2188.8} & \multicolumn{1}{c|}{9}    & \multicolumn{1}{c|}{1277.61} & \multicolumn{1}{c|}{9}     & \multicolumn{1}{c|}{1183.04}    & 0.0019                 &                      \\ 
\multicolumn{1}{|c|}{1920} & \multicolumn{1}{c|}{7680} & \multicolumn{1}{c|}{-}    & -    & \multicolumn{1}{c|}{-}    & \multicolumn{1}{c|}{-}    & \multicolumn{1}{c|}{-}    & \multicolumn{1}{c|}{-}    & \multicolumn{1}{c|}{10}     & \multicolumn{1}{c|}{13025.0}     &  \multicolumn{1}{c|}{10}     & \multicolumn{1}{c|}{12151.2}    &          0.0011             &                      \\ \cline{1-13}
\end{tabular}
}
\end{table}

\begin{table}[H]
\centering
\caption{\centering One-class {[}Non-Separable{]}: number of iterations for convergence (iter), root mean square error (RMSE), and CPU time (in seconds) for different spatial-temporal grid sizes $N_x \times N_t$. Comparison between sequential (LGMRES) and parallel results with different numbers of CPU cores. 
}
\label{tab_class1_petsc_nonseparable}
\resizebox{\textwidth}{!}{%
\begin{tabular}{cccc|ccccccccc|c}
\cline{5-13}
                           &                            &                           &      & \multicolumn{9}{c|}{parallel (CPU cores)}                                                                                                                                                                                                                        &                      \\ \cline{3-13}
                           & \multicolumn{1}{c|}{}      & \multicolumn{2}{c|}{LGMRES}       & \multicolumn{2}{c|}{1}                            & \multicolumn{2}{c|}{16}                          & \multicolumn{2}{c|}{56}                          & \multicolumn{2}{c|}{112}                        & \multirow{2}{*}{RMSE} &                      \\ \cline{1-12}
\multicolumn{1}{|c|}{$N_x$}   & \multicolumn{1}{c|}{$N_t$}    & \multicolumn{1}{c|}{iter} & time  & \multicolumn{1}{c|}{iter} & \multicolumn{1}{c|}{time}  & \multicolumn{1}{c|}{iter} & \multicolumn{1}{c|}{time}  & \multicolumn{1}{c|}{iter} & \multicolumn{1}{c|}{time}  & \multicolumn{1}{c|}{iter} & \multicolumn{1}{c|}{time} &                       &                      \\ \cline{1-13}
\multicolumn{1}{|c|}{30}   & \multicolumn{1}{c|}{120}   & \multicolumn{1}{c|}{4}    & 1.9892 & \multicolumn{1}{c|}{9}    & \multicolumn{1}{c|}{1.4183}  & \multicolumn{1}{c|}{9}    & \multicolumn{1}{c|}{0.3808}  & \multicolumn{1}{c|}{9}    & \multicolumn{1}{c|}{0.33785}  & \multicolumn{1}{c|}{9}     & \multicolumn{1}{c|}{0.36964}    & 0.0142                 &                      \\ 
\multicolumn{1}{|c|}{60}   & \multicolumn{1}{c|}{240}   & \multicolumn{1}{c|}{4}    & 8.6943 & \multicolumn{1}{c|}{11}   & \multicolumn{1}{c|}{8.4206}  & \multicolumn{1}{c|}{11}   & \multicolumn{1}{c|}{2.0786}  & \multicolumn{1}{c|}{11}   & \multicolumn{1}{c|}{1.66433}  & \multicolumn{1}{c|}{11}     & \multicolumn{1}{c|}{1.63943}    & 0.0086                 &                      \\ 
\multicolumn{1}{|c|}{120}  & \multicolumn{1}{c|}{480}   & \multicolumn{1}{c|}{5}    & 89.412 & \multicolumn{1}{c|}{15}   & \multicolumn{1}{c|}{56.405} & \multicolumn{1}{c|}{15}   & \multicolumn{1}{c|}{13.663} & \multicolumn{1}{c|}{15}   & \multicolumn{1}{c|}{11.0040} & \multicolumn{1}{c|}{15}     & \multicolumn{1}{c|}{10.2259}    & 0.0050                 &                      \\ 
\multicolumn{1}{|c|}{240}  & \multicolumn{1}{c|}{960}   & \multicolumn{1}{c|}{5}    & 5000.0 & \multicolumn{1}{c|}{17}   & \multicolumn{1}{c|}{350.53} & \multicolumn{1}{c|}{17}   & \multicolumn{1}{c|}{73.656} & \multicolumn{1}{c|}{17}   & \multicolumn{1}{c|}{59.1687} & \multicolumn{1}{c|}{17}     & \multicolumn{1}{c|}{52.8208}    & 0.0028                 &                      \\ 
\multicolumn{1}{|c|}{480}  & \multicolumn{1}{c|}{1920}  & \multicolumn{1}{c|}{-}    & -    & \multicolumn{1}{c|}{19}   & \multicolumn{1}{c|}{2417.9} & \multicolumn{1}{c|}{19}   & \multicolumn{1}{c|}{406.71} & \multicolumn{1}{c|}{19}   & \multicolumn{1}{c|}{317.431} & \multicolumn{1}{c|}{19}     & \multicolumn{1}{c|}{281.389}    & 0.0015                 &                      \\ 
\multicolumn{1}{|c|}{960}  & \multicolumn{1}{c|}{3840}  & \multicolumn{1}{c|}{-}    & -    & \multicolumn{1}{c|}{-}    & \multicolumn{1}{c|}{-}    & \multicolumn{1}{c|}{22}   & \multicolumn{1}{c|}{3445.0} & \multicolumn{1}{c|}{22}   & \multicolumn{1}{c|}{2021.92} & \multicolumn{1}{c|}{22}     & \multicolumn{1}{c|}{1775.82}    & 0.0008                &                      \\ 
\multicolumn{1}{|c|}{1920} & \multicolumn{1}{c|}{7680}  & \multicolumn{1}{c|}{-}    & -    & \multicolumn{1}{c|}{-}    & \multicolumn{1}{c|}{-}    & \multicolumn{1}{c|}{-}    & \multicolumn{1}{c|}{-}    & \multicolumn{1}{c|}{25}     & \multicolumn{1}{c|}{21228.1}     & \multicolumn{1}{c|}{25}     & \multicolumn{1}{c|}{20045.9}    &         0.0004              &                      \\ \cline{1-13}
\end{tabular}
}
\end{table}

We observe that for [LWR] in Table \ref{tab_class1_petsc_lwr}, by average $98 \%$ of the CPU time is saved by switching from sequential to parallel, $85 \%$ is saved for [Separable] in Table \ref{tab_class1_petsc_separable}, and $40 \%$ is saved for [Non-separable] in Table \ref{tab_class1_petsc_nonseparable}. Moving from $1$ to $16$ CPU cores, we save by an average of $71 \%$, $74 \%$, and $77 \%$ of time for [LWR], [Separable] and [Non-Separable], respectively. Adding more cores does not help much in terms of saved time since the percentage of improvement becomes smaller than $20 \%$. However, the PIST becomes greater in the finer grids. This latter requires more memory, which nevertheless justifies the need to use more CPU cores.\\

\amal{\noindent {\bf Future extension.} One significant challenge to be addressed in the future is related to the preconditioner since calculating the inverse of the approximate or exact Jacobian can be computationally expensive which heavily relies on the memory issues we faced. To address these limitations, future work should focus on developing more efficient and robust preconditioning strategies that could involve Diagonalization techniques, proposed in \cite{McDonald2018PreconditioningAI}.}


\section{Numerical Results}\label{section4}

In this section, we study the behavior of two classes of vehicles and their interaction using the game-theoretic framework approach proposed in section \ref{section2} and following the numerical methodologies proposed in section \ref{section3}.\\
\\
We consider two classes of vehicles on a single-lane road: class $1$ of light vehicles referred to as cars, and class $2$ of heavy vehicles referred to as trucks. Then it is natural to assume that $l_{1} < l_{2}$, $u_{max}^{1} > u_{max}^{2}$, \revision{$\rho^{1}_{jam} > \rho^{2}_{jam}$}.
\\
\\
To comply with the consistency relation \eqref{eq.density.lenth.relation}, the following inputs are selected:\\
For the class 1 of cars: $l_{1}=1$, $\rho^{1}_{jam}=1$, $u^{1}_{max}=1$, $L_1=1$. \\
For the class 2 of trucks: $l_{2}=2$, $\rho^{2}_{jam}=0.5$, $u^{2}_{max}=0.5$, $L_2=1$.\\
The initial densities for cars and trucks are denoted $\zeta_0^1$ and $\zeta_0^2$, respectively.

\subsection{Macroscopic simulations: NMFE}

\begin{figure}[H]
     \centering
     \includegraphics[width=\textwidth]{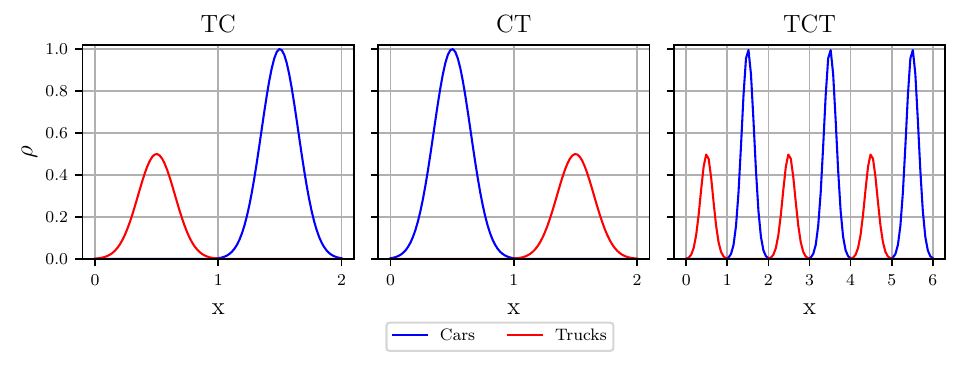}
    \caption{\centering Two-class initial density configurations: TC is a fully segregated configuration with \textcolor{blue}{cars} in front. CT is a fully segregated configuration with \textcolor{red}{trucks} in front. TCT is an interlaced configuration with alternating \textcolor{red}{trucks} and \textcolor{blue}{cars}.}
    \label{fig.2class_initial_config}
\end{figure}

We consider the time horizon $T=3$ and the following initial configurations \amal{(illustrated in Figure \ref{fig.2class_initial_config})}: 
\begin{itemize}
    \item Case 1: TC configuration. Fully segregated with cars in front with $\alpha_1=\alpha_2=1$. Then $L=L_1+L_2=2$, $\zeta_{0}^{1}(x)=\rho_{0,1}^{\left [ 1,2 \right ]}(x)$ and $\zeta_{0}^{2}(x)=\rho_{0,2}^{\left [ 0,1 \right ]}(x)$. 
    \item Case 2: CT configuration. Fully segregated with trucks in front with $\alpha_1=\alpha_2=1$. Then $L=L_1+L_2=2$, $\zeta_{0}^{1}(x)=\rho_{0,1}^{\left [ 0,1 \right ]}(x)$ and $\zeta_{0}^{2}(x)=\rho_{0,2}^{\left [ 1,2 \right ]}(x)$.
    \item Case 3: TCT configuration. Interlaced configuration with $\alpha_1=\alpha_2=3$. Then $L=3L_1+3L_2=6$, $\zeta_{0}^{1}(x)= \rho_{0,1}^{\left [ 1,2 \right ]}(x)+\rho_{0,1}^{\left [ 3,4 \right ]}(x)+\rho_{0,1}^{\left [ 5,6 \right ]}(x)$ and $\zeta_{0}^{2}(x)=\rho_{0,2}^{\left [0,1 \right ]}(x)+\rho_{0,2}^{\left [2,3 \right ]}(x)+\rho_{0,2}^{\left [4,5 \right ]}(x)$.
\end{itemize}
where the $j$-th class's initial density in $\left [ x_{1},x_{2} \right ]$ can be written as follow: 
\begin{equation}\label{eq.initial.density}
    \rho_{0,j}^{\left [ x_{1},x_{2} \right ]}(x) =\rho_{j}^{a}+(\rho_{j}^{b}-\rho_{j}^{a}) \cdot exp\left [ -\frac{\left ( x-\theta \right )^{2}}{2\gamma ^{2}} \right ] \cdot \mathbbm{1}_{\left [ x_{1},x_{2} \right ]}(x)   
\end{equation}
with $\rho_{1}^{a}=\rho_{2}^{a}=0.0$, $\rho_{1}^{b}=1.0$, $\rho_{2}^{b}=0.5$, $\gamma=0.15$ and $\theta =\frac{x_{1}+x_{2}}{2}$.\\

We refer to the three cases 1, 2, and 3, respectively, as TC, CT, and TCT in the presented figures.

\begin{figure}[H]
     \centering
     \includegraphics[width=\textwidth]{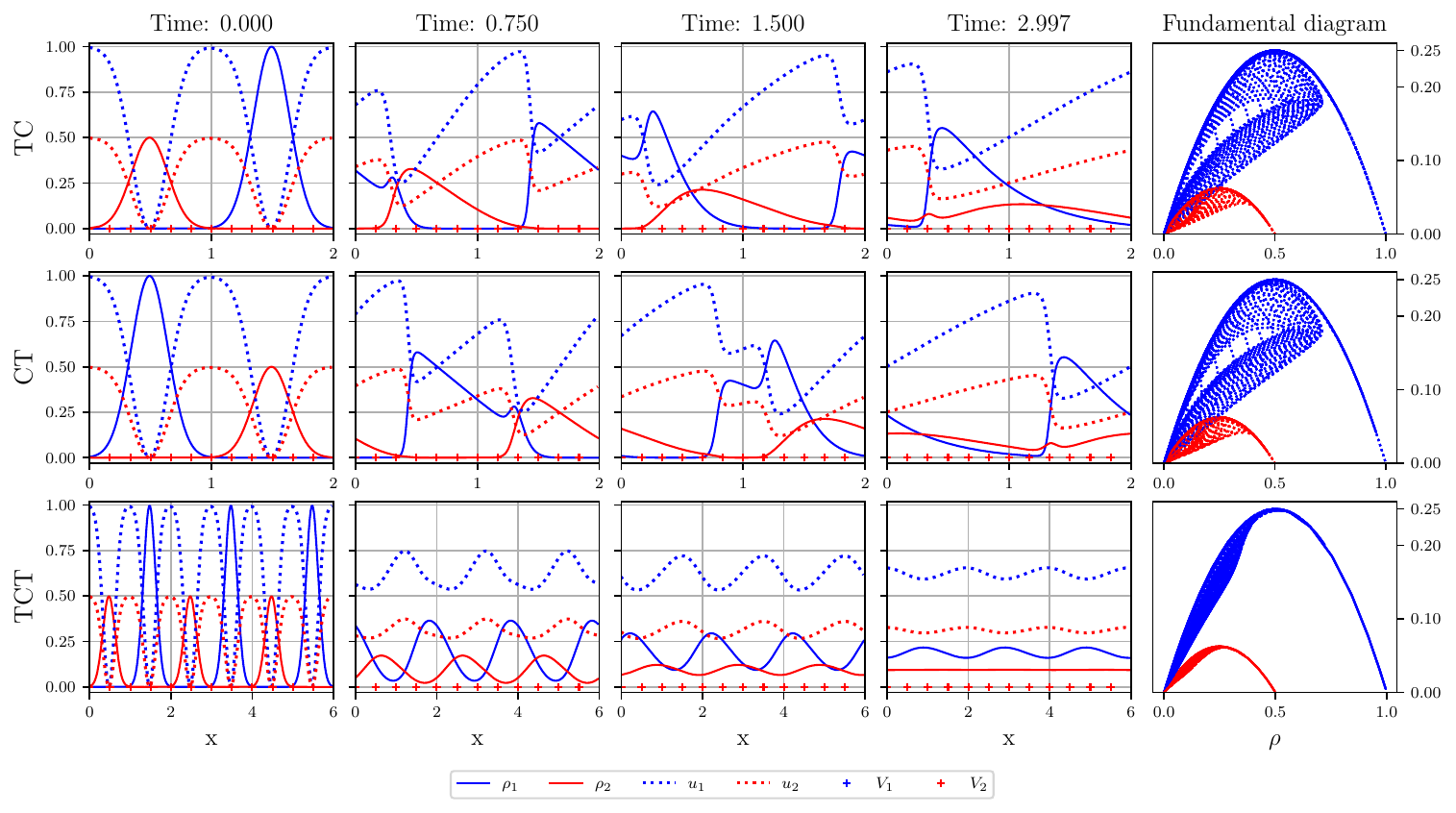}
      \caption{\centering Cost function [GLWR]: evolution of density $\rho_j$, speed $u_j$, and optimal cost $V_j$ for \textcolor{blue}{cars} and \textcolor{red}{trucks}, w.r.t. the position $x$. Each of the first four columns represents a time value, and the last column corresponds to the fundamental diagram. The rows represent the three cases TC, CT, and TCT.}
        \label{fig.macro_mfg_2class_LWR}
    
\end{figure}

\noindent {\bf Cost function 1: Generalized LWR.} \\
We consider a generalization of the Lighthill-Whitham-Richards, referred to in all what follows as [GLWR], for $n$-population model introduced in \cite{BenzoniGavage2003AnM}, where the speed of each class is influenced by the presence of the others, thus: 
$$u^{j}=U_{j}(\rho^{1},\rho^{2}), \text{for } j=1,2$$.\\
We consider the general Greenshields-type speed-density:
\begin{equation}\label{eq2.18}
    U_{j}(\rho^{1},\rho^{2})=u_{max}^{j}(1-s)
\end{equation}
with $s$ denoting the fraction of road occupancy ($0 \leq s \leq 1$), defined as follow:
\begin{equation}\label{eq.occupancy.1}
    s=\rho^{1}l_{1}+\rho^{2}l_{2}
\end{equation}
In other words, using \eqref{eq.density.lenth.relation}, and taking into consideration that $L_1=L_2=1$:
\begin{equation}\label{eq.occupancy.2}
    s=\frac{\rho^{1}}{\rho^{1}_{jam}}+\frac{\rho^{2}}{\rho^{2}_{jam}}
\end{equation}
Following the same lines as in the one-class case, we consider the cost function presented in \cite{huang_game-theoretic_2020}, which consists of keeping the actual speed of vehicles $u_{j}$ not too far from the desired one in \eqref{eq2.18}. Hence:
\begin{equation}\label{eq.generalized.lwr}
    f_{j}^{\text{GLWR}}(u^{j},\rho^{1},\rho^{2})=\frac{1}{2}(U_{j}(\rho^{1},\rho^{2})-u^{j})^{2}
\end{equation}\\
For [GLWR], the obtained density, speed, and optimal cost profiles for all scenarios, as shown in Figure \ref{fig.macro_mfg_2class_LWR} are qualitatively similar to those resulting from the single-class [LWR] model presented in \cite{huang_game-theoretic_2020} and reproduced in Figure \ref{fig.macro_mfg_1class}.

When cars and trucks are fully segregated (TC and CT) as shown in Figure \ref{fig.macro_mfg_2class_LWR}, we see that at some point they can manage to speed up and spread out their density around the mean. However, once one class reaches the other, we clearly observe a shock wave starting (Time: $0.75$), where the density before the shock point starts to increase. The shock is followed by a rarefaction along which the density after the shock point decreases. Because of the periodicity, the shock waves never disappear and are transported over time forming clusters of slow-moving vehicles, as shown at Time: $1.5$ and Time: $2.997$.\\
\\
The optimal costs $V_j$ remain null in Figure \ref{fig.macro_mfg_2class_LWR}, suggesting that for [GLWR] the optimal control precisely matches the desired value as required by the quadratic cost function in \eqref{eq.generalized.lwr}. As a consequence, the Greenshields-type law is attained by the mean-field macroscopic control.
\\
\\
In the interlaced configuration (TCT), as shown in Figure \ref{fig.macro_mfg_2class_LWR}, the traffic behavior is smoother. Indeed, no shock waves appear, and the density of both cars and trucks keeps spreading out over time. Although the correlation between the distribution of vehicles and the occurrence of shock waves is not clear, TCT experiment shows that a greater mixture can attenuate the formation of shock waves. This can be useful for traffic management, since more alternating distribution of vehicles is less likely to form clusters of slow-moving vehicles.\\
\\
We observe in the fundamental diagrams of Figure \ref{fig.macro_mfg_2class_LWR} that for all [GLWR] cases, three different regimes are distinguished: In the \textit{free-flow} regime, vehicles of both classes are in the free-flow ($\rho_{1} \leq 0.25$ and $\rho_{2} \leq 0.25$). In the \textit{semi-congested} regime, the trucks' class does not have enough space to maintain free flow, but the cars' class still does ($\rho_{1} \leq 0.5$ and $\rho_{2} \geq 0.25$). In the \textit{fully-congested} regime, both classes do not have sufficient space to maintain the free-flow speed ($\rho_{1} > 0.5$ and $\rho_{2} > 0.5$). In TC and CT configurations, most of the vehicles move to the semi-congested regime, which explains the formation of shock waves, unlike the TCT configuration, where most of the vehicles stay in the free-flow regime. But in the three configurations, only cars move to the fully-congested regime.\\
\\
Another important observed result for [GLWR] is that the fundamental relation satisfies all important requirements for the classical generalized LWR model in every class. For $j=1,2$:
\begin{itemize}
    \item Req1: Speed range $0 \leq v_{j} \leq v_{max}^{j}$  \label{enu1}
    \item Req2: Density range $0 \leq \rho _{j} \leq \rho _{jam}^{j}$  \label{enu2}
    \item Req3: No vehicle on the road $v_{j}=v_{max}^{j}$ when $\rho_{1}=\rho_{2}=0$ \label{enu3}
    \item Req4: Traffic jam $v_{j}=0$ when $\rho_{j}=\rho_{jam}^{j}$  \label{enu4}
    \item Req5: The fundamental diagram is strictly concave. \label{enu5}
\end{itemize}
Req3 accounts for the fact that here, the flow depends not only on density, as in one class case, but also on the composition of the traffic.\\

\noindent {\bf Cost function 2: Generalized Separable.} \\
We consider a generalization of the Separable cost function for one-class traffic proposed in \cite{huang_game-theoretic_2020}, referred to in all what follows as [GS], defined by:
\begin{equation}\label{eq.generalized.S.1}
    f_{j}^{\text{GS}}(u^{j},\rho^{1},\rho^{2})=\frac{1}{2}\left ( \frac{u^{j}}{u^{j}_{max}} \right )^{2}-\frac{u^{j}}{u^{j}_{max}}+ \frac{\rho}{\rho_{jam}}
\end{equation}
where,
$$\left\{\begin{matrix}
 \rho=\rho^{1}l_{1}+\rho^{2}l_{2} \\
 \rho_{jam}=\rho^{1}_{jam} l_{1}+\rho^{2}_{jam} l_{2}
\end{matrix}\right.$$
Using \eqref{eq.density.lenth.relation}, we get:
\begin{equation}\label{eq.rho}
  \left\{\begin{matrix}
 \rho=\frac{\rho^{1}}{\rho^{1}_{jam}}+\frac{\rho^{2}}{\rho^{2}_{jam}} \\
 \rho_{jam}=2
\end{matrix}\right.  
\end{equation}
Then \eqref{eq.generalized.S.1} becomes:
\begin{equation}\label{eq.generalized.S.2}
    f_{j}^{\text{GS}}(u^{j},\rho^{1},\rho^{2})=\frac{1}{2}\left ( \frac{u^{j}}{u^{j}_{max}} \right )^{2}-\frac{u^{j}}{u^{j}_{max}}+ \frac{1}{2} \left ( \frac{\rho^{1}}{\rho^{1}_{jam}}+\frac{\rho^{2}}{\rho^{2}_{jam}} \right )
\end{equation}
\begin{figure}[H]
     \centering
     \includegraphics[width=\textwidth]{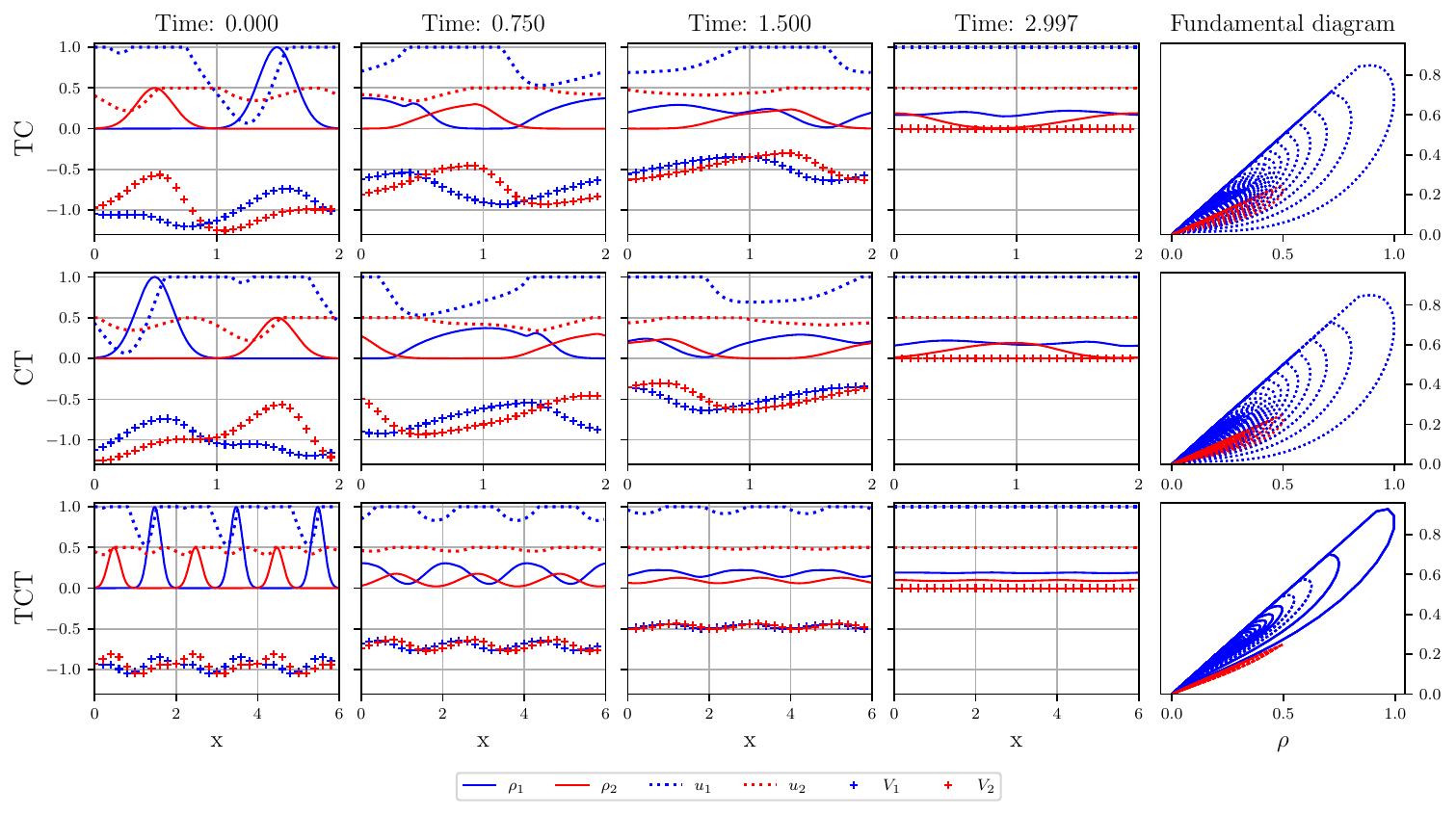}
    \caption{\centering Cost function [GS]: evolution of density $\rho_j$, speed $u_j$, and optimal cost $V_j$ for \textcolor{blue}{cars} and \textcolor{red}{trucks}, w.r.t. the position $x$. Each of the first four columns represents a time value, and the last column corresponds to the fundamental diagram. The rows represent the three cases TC, CT, and TCT.}
    \label{fig.macro_mfg_2class_sep}
\end{figure}

\noindent {\bf Cost function 3: Generalized Non-Separable.} \\
We consider now a generalization of the Non-Separable cost function for one-class traffic proposed in \cite{huang_game-theoretic_2020}, referred to in all what follows as [GNS] defined by:
\begin{equation}\label{eq.generalized.NS.1}
    f_{j}^{\text{GNS}}(u^{j},\rho^{1},\rho^{2})=\frac{1}{2}\left ( \frac{u^{j}}{u^{j}_{max}} \right )^{2}-\frac{u^{j}}{u^{j}_{max}}+\frac{u^{j}}{u^{j}_{max}} \frac{\rho}{\rho_{jam}}
\end{equation}
Using \eqref{eq.rho}, the cost \eqref{eq.generalized.NS.1} becomes:
\begin{equation}\label{eq.generalized.NS.2}
    f_{j}^{\text{GNS}}(u^{j},\rho^{1},\rho^{2})=\frac{1}{2}\left ( \frac{u^{j}}{u^{j}_{max}} \right )^{2}-\frac{u^{j}}{u^{j}_{max}}+\frac{1}{2} \frac{u^{j}}{u^{j}_{max}} \left ( \frac{\rho^{1}}{\rho^{1}_{jam}}+\frac{\rho^{2}}{\rho^{2}_{jam}} \right )
\end{equation}
\begin{figure}[H]
     \centering
         \includegraphics[width=\textwidth]{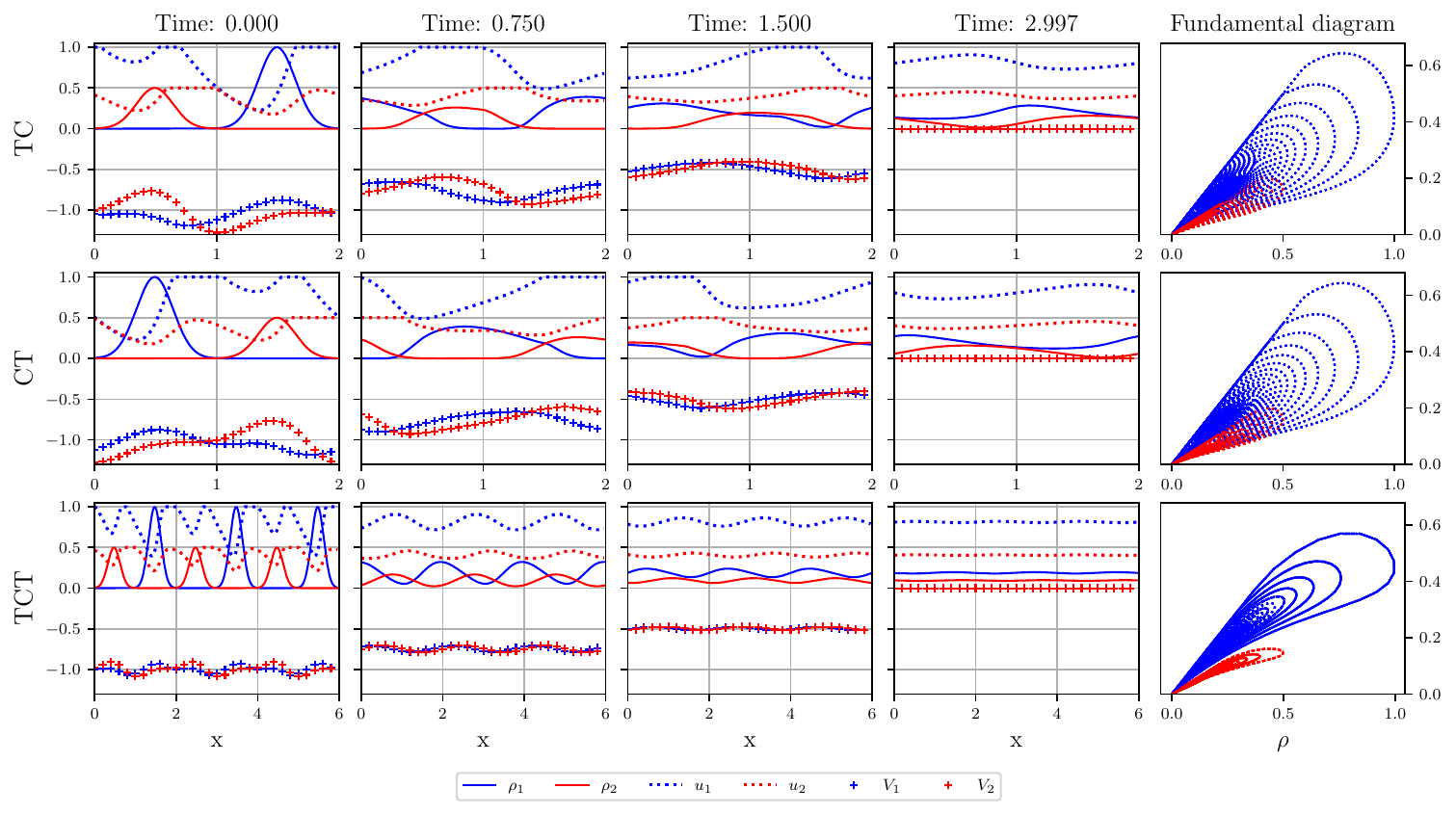}
        \caption{\centering Cost function [GNS]: evolution of density $\rho_j$, speed $u_j$, and optimal cost $V_j$ for \textcolor{blue}{cars} and \textcolor{red}{trucks}, w.r.t. the position $x$. Each of the first four columns represents a time value, and the last column corresponds to the fundamental diagram. The rows represent the three cases TC, CT, and TCT.}
    \label{fig.macro_mfg_2class_nonsep}
\end{figure}
For [GS] in Figure \ref{fig.macro_mfg_2class_sep}, and [GNS] in Figure \ref{fig.macro_mfg_2class_nonsep}, there is no shock formation, thanks to the third term in the cost functions \eqref{eq.generalized.S.2} and \eqref{eq.generalized.NS.2} that takes into account the level of density area as well as the composition of traffic. We observe the anticipation behavior in which the vehicles favor driving at a lower speed than the desired one, even in a low-density area, if there is a traffic jam ahead, contrary to [GLWR], where vehicles try to get as close as possible to the desired speed without anticipating the capacity of the road. As a result, the vehicles in [GLWR] are more sensitive to the changes in density, which clearly appears in the fundamental diagram of Figure \ref{fig.macro_mfg_2class_LWR} where we have steeper curves indicating that their flows are strongly impacted by density evolution. This information can help to identify the vehicles' behaviors that are contributing most to congestion and to design targeted strategies that can be useful for congestion mitigation.\\
\\
Contrary to [GLWR], the optimal cost for both [GS] in Figure \ref{fig.macro_mfg_2class_sep} and [GNS] in Figure \ref{fig.macro_mfg_2class_nonsep} does not vanish and varies within the range $\left [ -1.5 , 0.0 \right ]$, indicating that the system is trying to balance the different terms of $f_j$ in order to find the optimal velocity control. The fact that the optimal cost takes on negative values implies that the negative term $-\frac{u^{j}}{u^{j}_{max}}$ is meaningful. This term penalizes the linear deviation from the maximum speed.\\
\\
In both [GS] and [GNS], there is no fully-congested regime. Instead, a synchronized free-flow regime appears at Time: $2.997$, where vehicles of the same class move at a constant maximum speed and with uniform distribution, resulting in a smooth flow of traffic. Hence, the special closed shape of the fundamental diagrams in Figure \ref{fig.macro_mfg_2class_sep}, and Figure \ref{fig.macro_mfg_2class_nonsep}.\\

\subsection{Microscopic discrete controls}
In this section, we vary the number of vehicles in each class $(N_j)_{j=1,2}$ such that $N_j=\alpha_j n$, with $n \in \left\{ 20, 40, 60, 80, 100\right\}$. Let $\alpha_1=\alpha_2=1$ for TC and CT, and $\alpha_1=\alpha_2=3$ for TCT.\\
We assume $L_j=n$, then $\rho_{jam}^{j}$ should fulfill \eqref{eq.density.lenth.relation}, and $L=N_1 + N_2$.\\
To reach the stable state at $T=3$, we consider the $j$-th class free flow speed $u_{max}^{j}\times L_j$.\\
\\
Following the micro-macro approach introduced in Section \ref{sec2.2}, we first compute the NMFG-constructed controls from the macroscopic NMFG solutions by solving \eqref{eq.nmfg.constructed.control}, using Algorithm \ref{alg:NMFE-constructed controls}.
\begin{small}
\begin{algorithm}[H]
\caption{NMFE-constructed controls}
\label{alg:NMFE-constructed controls}
\begin{algorithmic}[1]
    \State Given $N_{1}$ the number of cars and $N_{2}$ the number of trucks. $N=N_1 + N_2$ the total number of vehicles. $(u_{1}^{*},u_{2}^{*})$ and $(\rho_{1}^{*},\rho_{2}^{*})$ the NMFE solutions.
    \State Generate from the same initial densities used in macro-level $(\rho_{1}^{0},\rho_{2}^{0})$ : $N_{1}$ random samples for cars denoted $x_{1,0},...,x_{N_{1},0}$, and $N_{2}$ random samples for trucks denoted $x_{N_{1}+1,0},...,x_{N,0}$. 
    \State Compute the NMFE-constructed controls {$\left\{\hat{v}_{i_{j}}(t) \right\}_{i_j \in \mathcal{I}_j,j=1,2}$} using \eqref{eq.nmfg.constructed.control}.
    \State Compute the best response strategies {$\left\{\bar{v}_{i_{j}}(t) \right\}_{i_j \in \mathcal{I}_j,j=1,2}$} using algorithm \ref{alg:Best_response_strategie}.
    \State Compute the accuracy associated {$\left\{\hat{\varepsilon}_{i_{j}} \right\}_{i_j \in \mathcal{I}_j,j=1,2}$} using \eqref{eq.epsilon.lower.bound}.
    \State Compute the maximum and mean accuracy :
    $$\left\{\begin{matrix}
    \text{MaxRA}=\frac{\underset{j=1,2}{max}\left\{ \underset{i_j \in \mathcal{I}_j}{max}\left|\hat{\varepsilon}_{i_{j}} \right|\right\}}{\underset{j=1,2}{max}\left\{\underset{i_j\in \mathcal{I}_j}{max}\left|J^N_{i_{j}}(\hat{v}_{i_{j}},\hat{v}_{-i_{j}},\hat{v}^{-j}) \right|\right\}}\\
    \text{MeanRA}=\frac{\sum_{j=1,2}^{}\sum_{i_j \in \mathcal{I}_j}^{} \left|\hat{\varepsilon}_{i_{j}} \right|}{\sum_{j=1,2}^{}\sum_{i_j\in \mathcal{I}_j}^{} \left|J^N_{i_{j}}(\hat{v}_{i_{j}},\hat{v}_{-i_{j}},\hat{v}^{-j})\right|}
    \end{matrix}\right.$$
\end{algorithmic}
\end{algorithm}
\end{small}

We then compare them with the best response strategies computed by solving equations \eqref{eq.best.response.strategy}-\eqref{eq.optimal.dynamic}, using Algorithm \ref{alg:Best_response_strategie}.
\begin{algorithm}[H]
\caption{Best response strategy}
\label{alg:Best_response_strategie}
\begin{algorithmic}[1]
    \For{$j: 1,2$}
    \For{$i_j \in \mathcal{I}_j$}
    \State Initialize $k=0$, and choose a fixed line search step $\tau > 0$ and a convergence tolerance $\epsilon<<1$. 
    \State Choose $v_{i_j}^{(0)}(t)\in \left [ 0,1 \right ]$ (e.g, $v_{i_j}^{(0)}(t)=\frac{\hat{v}_{i_j}(t)}{u^j_{max} \times L_j}$).
    \For{$m \in \mathcal{I}_1 \cup \mathcal{I}_2$}
    \State Solve the state equation (Explicit forward Euler) :
            $$\left\{\begin{matrix}
        \dot{x}_m^{(k)}(t)=\left\{\begin{matrix}
        \hat{v}_m(t)/(u^j_{max} \times L_j) & \text{if } m \neq i_j \\
        v_m^{(k)}(t)/(u^j_{max} \times L_j) & \text{if } m = i_j  \\
        \end{matrix}\right.\\
        x_m^{(k)}(0)=x_{m,0}/L
        \end{matrix}\right.$$
    \EndFor
    \State Solve the backward adjoint state equation (Explicit forward Euler) : 
            $$\left\{\begin{matrix}
        \dot{P}_{i_j}^{(k)}(t)= -\nabla_x f_j\left ( v_{i_j}^{(k)}(t),\rho^1\left ( t,x_{i_j}^{(k)}(t) \right ),\rho^2\left ( t,x_{i_j}^{(k)}(t) \right ) \right )
        \\
        P_{i_j}^{(k)}(T)=0
        \end{matrix}\right.$$
    \State Compute the gradient :
    $$g_{i_j}^{(k)}(t)=P_{i_j}^{(k)}(t) + \nabla_v f_j\left ( v_{i_j}^{(k)}(t),\rho^1\left ( t,x_{i_j}^{(k)}(t) \right ),\rho^2\left ( t,x_{i_j}^{(k)}(t) \right ) \right )  $$
    \If{$\| g_{i_j}^{(k)}(t) \| \leq \epsilon$} 
        \State $v_{i_j}^{(k+1)}(t)=v_{i_j}^{(k)}(t)$
    \Else
        \State $v_{i_j}^{(k+1)}(t)=v_{i_j}^{(k)}(t)-\tau g_{i_j}^{(k)}(t)$
    \EndIf
    \If{$v_{i_j}^{(k+1)}(t) \approx v_{i_j}^{(k)}(t)$} 
        \State $\bar{v}_{i_j}(t)=v_{i_j}^{(k+1)}(t) \times u^j_{max} \times L_j$
    \Else
        \State $k=k+1$. 
        \State Repeat from 4:
    \EndIf
    \EndFor
    \EndFor

\end{algorithmic}
\end{algorithm}
\begin{figure}[H]
     \centering
         \includegraphics[width=\textwidth]{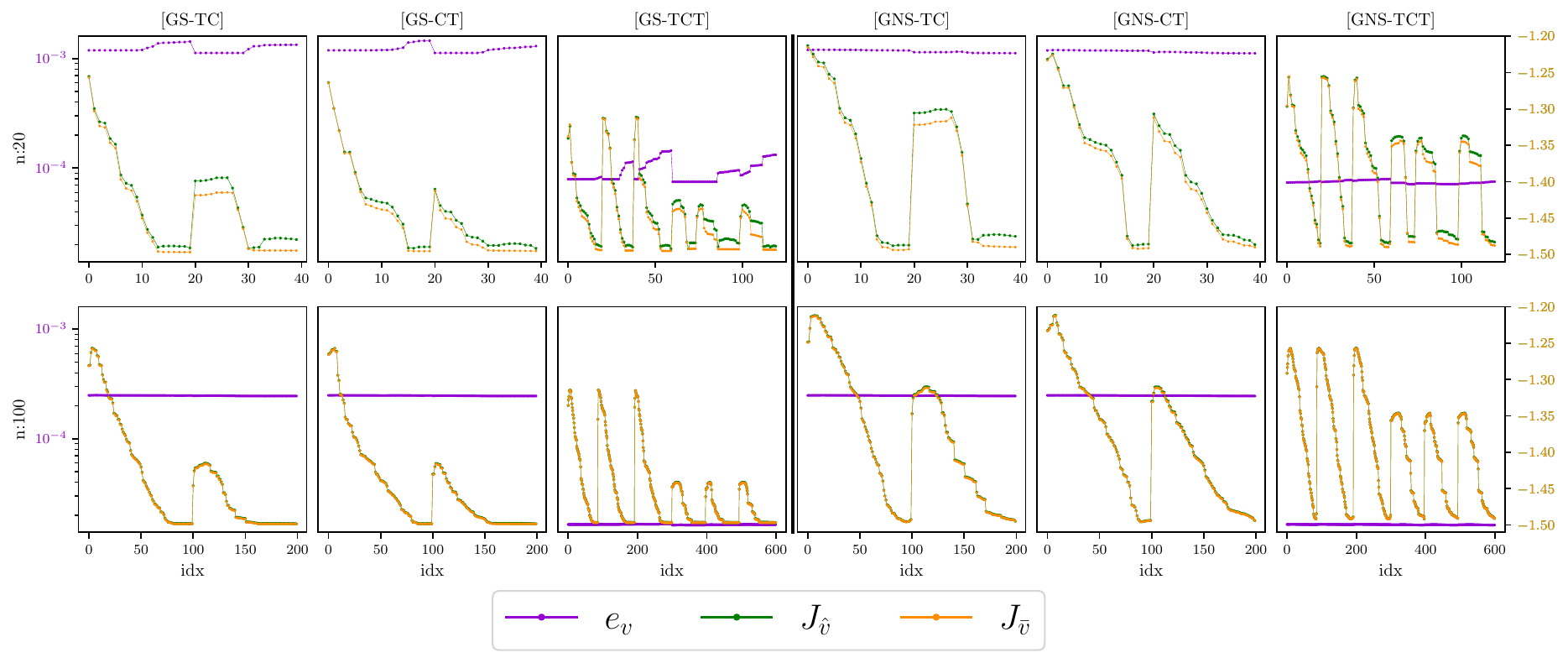}
        \caption{\centering The columns correspond to the TC, CT, and TCT configurations, considering [GS] (left) and [GNS] (right). The rows are for two different values of $n=20,100$. The total number of vehicles $N$ for each column is such that $N=N_1+N_2$, where $N_1=N_2=\alpha n$, with $\alpha=1$ for TC and CT and $\alpha=3$ for TCT. Each sub-figure presents three curves: the $L_{\infty}$ norm \textcolor{darkpurple}{$e_v=\left\|\hat{v}-\Bar{v}\right\|_{\infty}$}, the cost for NMFE-constructed controls \textcolor{darkgreen}{$J_{\hat{v}}=J_{i_{j}}(\hat{v}_{i_{j}},\hat{v}_{-i_{j}},\hat{v}^{-j})$}, and the cost for the best response strategies \textcolor{darkorange}{$J_{\bar{v}}=J_{i_{j}}(\Bar{v}_{i_{j}},\hat{v}_{-i_{j}},\hat{v}^{-j})$}, all w.r.t. vehicle's index (idx).}
    \label{fig.micro_2class_micro_plots}
\end{figure}

\begin{figure}[H]
     \centering
         \includegraphics[width=\textwidth]{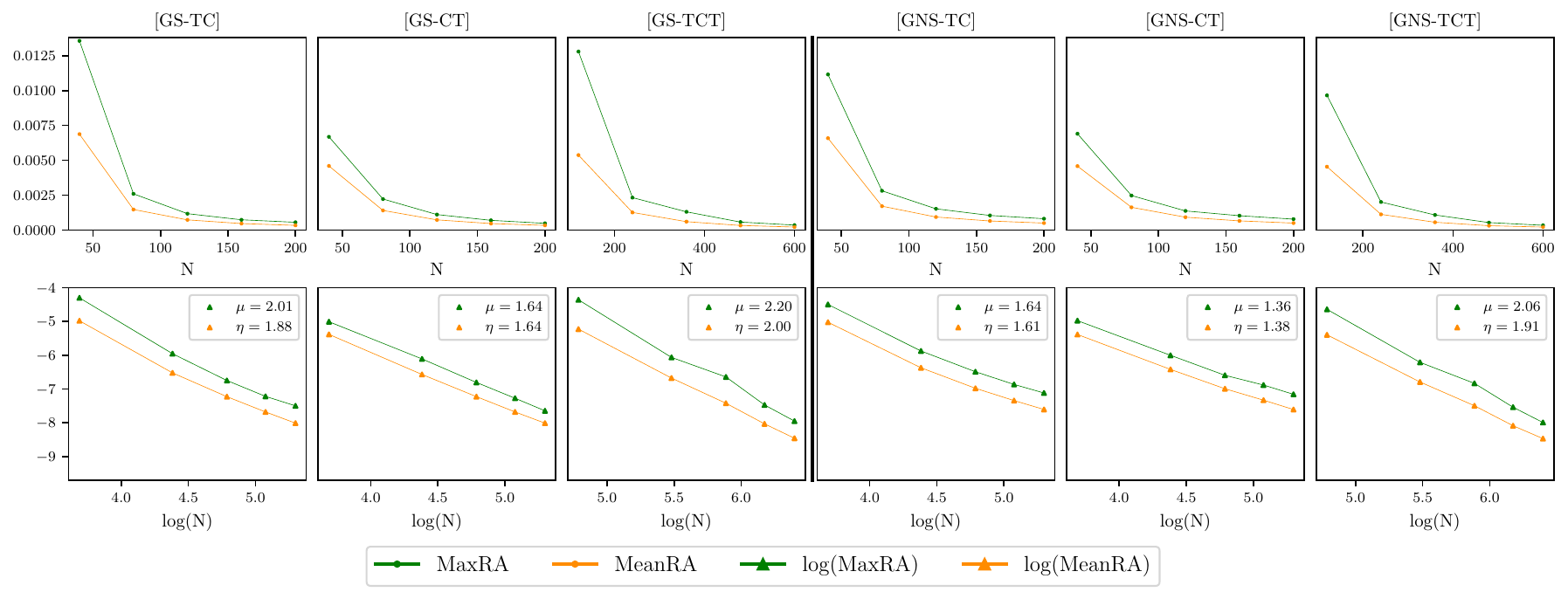}
        \caption{\centering The columns correspond to the TC, CT, and TCT configurations, considering [GS] (left) and [GNS] (right). The first row represents the mean accuracy (\textcolor{darkorange}{MeanRA}) and maximum accuracy (\textcolor{darkgreen}{MaxRA}) for different total number of vehicles $N$. The second row shows their respective log-log plots, with $-\mu$ and $-\eta$ being the slopes.}
    \label{fig.micro_2class_RA}
\end{figure}

From Figure \ref{fig.micro_2class_micro_plots}, we notice the following results: the values of cost function when vehicles move with their NMFE-constructed controls  $J_{\hat{v}}=J_{i_{j}}(\hat{v}_{i_{j}},\hat{v}_{-i_{j}},\hat{v}^{-j})$ are closer to those when vehicles move with their best response strategies  $J_{\bar{v}}=J_{i_{j}}(\Bar{v}_{i_{j}},\hat{v}_{-i_{j}},\hat{v}^{-j})$ which proves that the \revision{micro-macro} approach provides a good approximation of the vehicles control and shows the relevance of the algorithms proposed. We observe that the $L_{\infty}$ norm of the difference between the NMFE-constructed controls $\hat{v}$ and the best response strategies $\Bar{v}$ denoted by $e_v=\left\|\hat{v}-\Bar{v}\right\|_{\infty}$ are small and decreases every time the total number of vehicles $N$ increases, which is consistent with the main idea of the mean-field theory. In order to keep the paper concise, we showcase in Figure \ref{fig.micro_2class_micro_plots} the results for only two different $N$ values since the convergence is later outlined in Figure \ref{fig.micro_2class_RA}. $N$ is such that $N=N_1+N_2$, where $N_1=N_2=\alpha n$, where $n=20,100$, with $\alpha=1$ for TC and CT and $\alpha=3$ for TCT. 

\amal{The plots in the first row of Figure \ref{fig.micro_2class_RA} clearly show that the meanRA and maxRA are close to zero. Moreover, an increase in the number of vehicles $N$ leads to a decrease in meanRA and maxRA, implying that the accuracy $\hat{\varepsilon}_{i_j}$ is close to zero for all $i_j$, which is also clear for every vehicle in Figure \ref{fig.micro_2class_micro_plots}, meaning the NMFE-constructed control $\hat{v}_{i_j}$ are becoming increasingly close to the best strategy $\Bar{v}_{i_{j}}$ as $N$ increases. }
Moving to the second row, where log-log plots are displayed for each configuration to analyze the relationship between the relative accuracies and $N$. We clearly observe that MaxRA and MeanRA behave as $N^{-\mu}$ and $N^{-\eta}$ respectively. The experimental results found for $\mu$ and $\eta$ are larger than $\frac{1}{4}$, the latter is what is expected by the theory, see \cite{carmona2018probabilistic} section 6.1 Theorem 6.7. The same behavior appears in the case of one-class in Figure \ref{fig.micro_1class_RA}.
\begin{figure}[!ht]
     \centering
         \includegraphics[width=\textwidth]{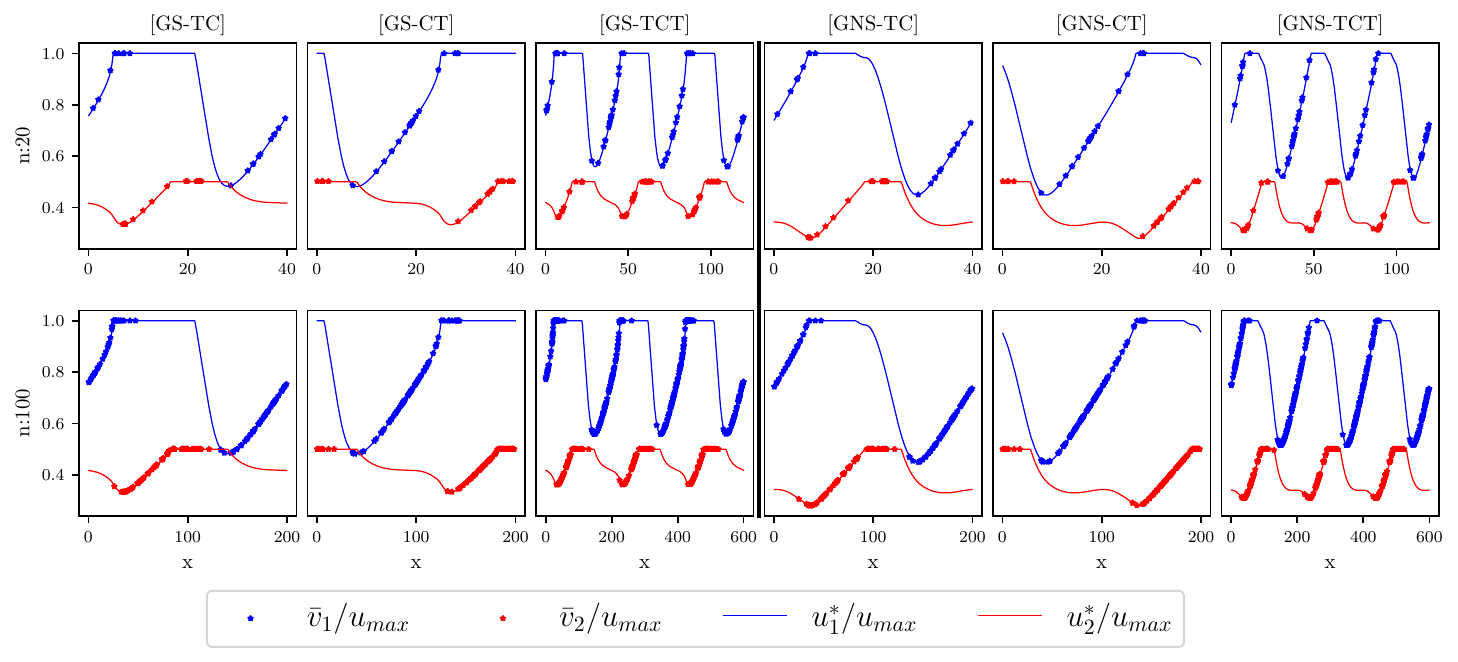}
        \caption{\centering The columns correspond to the TC, CT, and TCT configurations, considering [GS] (left) and [GNS] (right). The rows are for two different values of $n=20,100$. The total number of vehicles $N$ for each column is such that $N=N_1+N_2$, where $N_1=N_2=\alpha n$, with $\alpha=1$ for TC and CT and $\alpha=3$ for TCT. Each sub-figure presents: the relative macro NMFE $u^*_j/u_{max}$ and the relative best controls $\bar{v}_j/u_{max}$ at $t=0.6$, for \textcolor{blue}{cars} (j=1) and \textcolor{red}{trucks} (j=2).}
    \label{fig.micro_2class_control_t0.6}
\end{figure}

Figure \ref{fig.micro_2class_control_t0.6} shows the anticipation behavior of vehicles on the microscopic scale for [GS] and [GNS], which is consistent with what is described on the macroscopic scale in Figure \ref{fig.macro_mfg_2class_sep} and Figure \ref{fig.macro_mfg_2class_nonsep}. Again, we showcase in Figure \ref{fig.micro_2class_control_t0.6} the results for only two different $n$ values since the behavior is the same for the other values between $20$ and $100$.
\newpage
\section{Conclusion}

In this paper, we have explored a multi-class traffic model and examined the computational feasibility of mean-field games (MFG) in obtaining approximate Nash equilibria for traffic flow games involving a large number of players. We introduced a two-class traffic framework, building upon the ideas presented in the single-class approach by \cite{huang_game-theoretic_2020}. To facilitate our analysis, we employed various numerical techniques, including high-performance computing and regularization of LGMRES solvers. By utilizing these tools, we conducted simulations at significantly larger spatial and temporal scales.

We led extensive numerical experiments considering three different scenarios involving cars and trucks, as well as three different cost functionals. Our results primarily focused on the dynamics of autonomous vehicles (AVs) in traffic, which support the effectiveness of the approach proposed by \cite{huang_game-theoretic_2020}.

Moreover, we conducted original comparisons between macroscopic Nash mean-field speeds and their microscopic counterparts. These comparisons allowed us to computationally validate the $\epsilon-$Nash approximation, demonstrating a slightly improved convergence rate compared to theoretical expectations.

Future directions encompass second order traffic models, the multi-lane case, particularly prone to non-cooperative game considerations, and addressing some theoretical issues. \\

Actually, first-order models have severe limitations in describing important traffic features, such as light traffic with slow drivers and stop-and-start transitions. The first second-order model was introduced in \cite{payne1971model}, which at its turn was criticized in \cite{daganzo_requiem_1995} for it violates critical traffic principles such as that cars react to front stimuli (and not to cars behind). The model was "resurrected" later in \cite{aw_resurrection_2000,zhang_non-equilibrium_2002}. However, the game version of the model has not yet been explored.
 Moreover, this work may be extended to handle the case where the traffic flow is on a multi-lane road, for which non-game model versions are proposed on both microscopic and macroscopic scales. See for example \cite{YangQi1996,nagalur2019first}.  
 Finally, numerous theoretical issues arising from the two-class modeling, such as existence of solutions to the generalized limit systems (e.g. generalized LWR and generalized separable systems) or convergence properties of the discretized nonlinear systems.\\
\\
{\bf Declarations. }This work was supported in part by using the resources of Mohammed VI Polytechnic University through Simlab \& Toubkal the HPC \& IA platforms.\\
The authors have no financial or proprietary interests in any material discussed in this article.\\
The authors have no competing interests.\\
{\bf Data availability statement.} The datasets generated during and/or analysed during the current study are available from the corresponding author on reasonable request.\\
The datasets generated during and/or analysed during the current study are available in the GitHub link \url{https://github.com/amalmach/MFG_Traffic}, under MFG\_Traffic/ repository.

\printbibliography

@article{festa2018mean,
  title={A mean field game approach for multi-lane traffic management},
  author={Festa, Adriano and G{\"o}ttlich, Simone},
  journal={IFAC-PapersOnLine},
  volume={51},
  number={32},
  pages={793--798},
  year={2018},
  publisher={Elsevier}
}

@article{cirant2015multi,
  title={Multi-population mean field games systems with Neumann boundary conditions},
  author={Cirant, Marco},
  journal={Journal de Math{\'e}matiques Pures et Appliqu{\'e}es},
  volume={103},
  number={5},
  pages={1294--1315},
  year={2015},
  publisher={Elsevier}
}

@article{rahman2023impacts,
  title={Impacts of connected and autonomous vehicles on urban transportation and environment: A comprehensive review},
  author={Rahman, Md Mokhlesur and Thill, Jean-Claude},
  journal={Sustainable Cities and Society},
  volume={96},
  pages={104649},
  year={2023},
  publisher={Elsevier}
}

@article{payne1971model,
  title={Model of freeway traffic and control},
  author={Payne, Harold J},
  journal={Mathematical Model of Public System},
  pages={51--61},
  year={1971}
}

@article{YangQi1996,
  title={ A microscopic traffic simulator for evaluation of dynamic
traffic management systems},
  author={Yang, Qi and Koutsopoulos H. N.},
  journal={ Transportation Research Part C: Emerging Technology},
  volume={4(3)},
  pages={113-129},
  year={1996}
}

@article{lasry_mean_2007,
	title = {{Mean field games}},
	volume = {2},
	pages = {229--260},
	journaltitle = {{Japanese Journal of Mathematics}},
	author = {Lasry, Jean-Michel and Lions, Pierre-Louis},
	date = {2007-03-01},
	number = {{1}}
}

@article{huang_large_2006,
	title = {Large population stochastic dynamic games: Closed-loop {McKean}-Vlasov systems and the Nash certainty equivalence principle},
	volume = {6},
	shorttitle = {Large population stochastic dynamic games},
	journal = {{Communications in Information \& Systems}},
    number = {3},
    pages = {221 -- 252},
	author = {Minyi Huang and Roland P. Malhamé and Peter E. Caines},
	date = {2006-01-01},
}

@article{bensoussan_mean_2018,
	title = {{Mean Field Control and Mean Field Game Models with Several Populations}},
	journal = {{Minimax Theory and its Applications}},
	author = {Bensoussan, Alain and Huang, Tao and Lauriere, Mathieu},
	date = {2018-10-29},
	volume = {{3}},
	number = {{2, SI}},
	pages = {{173-209}}
}

@article{daganzo_requiem_1995,
	title = {Requiem for second-order fluid approximations of traffic flow},
	volume = {29},
	issn = {0191-2615},
	url = {https://www.sciencedirect.com/science/article/pii/019126159500007Z},
	doi = {10.1016/0191-2615(95)00007-Z},
	pages = {277--286},
	number = {4},
	journaltitle = {Transportation Research Part B: Methodological},
	shortjournal = {Transportation Research Part B: Methodological},
	author = {Daganzo, Carlos F.},
	date = {1995-08-01},
	langid = {english},
}

@article{aw_resurrection_2000,
	title = {Resurrection of "Second Order" Models of Traffic Flow},
	volume = {60},
	issn = {0036-1399},
	doi = {10.1137/S0036139997332099},
	pages = {916--938},
	number = {3},
	journaltitle = {{SIAM} Journal on Applied Mathematics},
	shortjournal = {{SIAM} J. Appl. Math.},
	author = {Aw, A. and Rascle, M.},
	date = {2000-01-01},
}

@article{zhang_non-equilibrium_2002,
	title = {A non-equilibrium traffic model devoid of gas-like behavior},
	volume = {36},
	issn = {0191-2615},
	doi = {10.1016/S0191-2615(00)00050-3},
	pages = {275--290},
	number = {3},
	journaltitle = {Transportation Research Part B: Methodological},
	shortjournal = {Transportation Research Part B: Methodological},
	author = {Zhang, H. M.},
	date = {2002-03-01},
	langid = {english},
}

@article{lauriere_numerical_2021,
	title = {Numerical Methods for Mean Field Games and Mean Field Type Control},
	journaltitle = {{arXiv}:2106.06231 [cs, math]},
	author = {Lauriere, Mathieu},
	date = {2021-06-11},
}

@article{huang_game-theoretic_2020,
	title = {A Game-Theoretic Framework for Autonomous Vehicles Velocity Control: Bridging Microscopic Differential Games and Macroscopic Mean Field Games},
	author = {Huang, Kuang and Di, Xuan and Du, Qiang and Chen, Xi},
	volume = {25},
	pages = {4869--4903},
	number = {12},
	journal = {{Discrete and Continuous Dynamical Systems-Series B}},
	date = {2020}
}

@article{knoll_jacobian-free_2004,
	title = {{Jacobian-free Newton-Krylov methods: a survey of approaches and applications}},
	volume = {193},
	shorttitle = {Jacobian-free Newton-Krylov Methods},
	number = {{2}},
	pages = {357--397},
	journaltitle = {{Journal of Computational Physics}},
	author = {Knoll, Dana and Keyes, David},
	date = {2004-01-20}}

@article{BenzoniGavage2003AnM,
    author = {Sylvie Benzoni-Gavage and Rinaldo M. Colombo},
    title = {{An n-populations model for traffic flow}},
    journal = {{Europ. J. Appl.Math}},
    year = 2003,
    volume = {{14}},
    number = {{5}},
    pages = {{587-612}},
    month = OCT  }

@article{ WOS:000231357700001,
    author = {Lewis, RM and Nash, SG},
    title = {{Model problems for the multigrid optimization of systems governed by
       differential equations}},
    journal = {{SIAM Journal on Scientific Computing}},
    year = 2005,
    volume = {{26}},
    number = {{6}},
    pages = {{1811-1837}}
}

@article{ WOS:000230502200005,
    author = {Baker, AH and Jessup, ER and Manteuffel, T},
    title = {{A technique for accelerating the convergence of restarted GMRES}},
    journal = {{SIAM Journal on Matrix Analysis and Applications}},
    year = 2005,
    volume = {{26}},
    number = {{4}},
    pages = {{962-984}}
}

@article{Achdou2020MeanFG,
  title={Mean field games and applications: Numerical aspects},
  author={Achdou, Yves and Lauri{\`e}re, Mathieu},
  journal={Mean field games},
  pages={249--307},
  year={2020},
  publisher={Springer}
}

@book{achdou2013hamiltonbook,
  title={Hamilton-Jacobi Equations: Approximations, Numerical Analysis and Applications: Cetraro, Italy 2011, Editors: Paola Loreti, Nicoletta Anna Tchou},
  author={Achdou, Y. and Barles, G. and Ishii, H. and Litvinov, G.L. and Loreti, P. and Tchou, N.},
  isbn={9783642364334},
  lccn={2013937326},
  series={Lecture Notes in Mathematics},
  url={https://books.google.co.ma/books?id=MiS6BQAAQBAJ},
  year={2013},
  publisher={Springer Berlin Heidelberg}
}

@article{Dalcin2021mpi4pySU,
  title={mpi4py: Status Update After 12 Years of Development},
  author={Lisandro Dalcin and Yao-Lung L. Fang},
  journal={Computing in Science \& Engineering},
  year={2021},
  volume={23},
  pages={47-54}
}

@misc{petsc-web-page,
  author = {Satish Balay and Shrirang Abhyankar and Mark~F. Adams and Steven Benson and Jed Brown
    and Peter Brune and Kris Buschelman and Emil~M. Constantinescu and Lisandro Dalcin and Alp Dener
    and Victor Eijkhout and Jacob Faibussowitsch and William~D. Gropp and V\'{a}clav Hapla and Tobin Isaac and Pierre Jolivet
    and Dmitry Karpeev and Dinesh Kaushik and Matthew~G. Knepley and Fande Kong and Scott Kruger
    and Dave~A. May and Lois Curfman McInnes and Richard Tran Mills and Lawrence Mitchell and Todd Munson
    and Jose~E. Roman and Karl Rupp and Patrick Sanan and Jason Sarich and Barry~F. Smith
    and Stefano Zampini and Hong Zhang and Hong Zhang and Junchao Zhang},
  title        = {{PETS}c {W}eb page},
  url          = {https://petsc.org/},
  howpublished = {\url{https://petsc.org/}},
  year         = {2022},
}

@article{bourne2023pyccel,
   title={Pyccel: a Python-to-X transpiler for scientific high-performance computing},
   author={Bourne, Emily and G{\"u}{\c{c}}l{\"u}, Yaman and Hadjout, Said and Ratnani, Ahmed},
   journal={Journal of Open Source Software},
   volume={8},
   number={83},
   pages={4991},
   year={2023}
 }

@article{abhyankar2018petsc,
   title={PETSc/TS: A modern scalable ODE/DAE solver library},
   author={Abhyankar, Shrirang and Brown, Jed and Constantinescu, Emil M and Ghosh, Debojyoti and Smith, Barry F and Zhang, Hong},
   journal={arXiv preprint arXiv:1806.01437},
   year={2018}
 }

@techreport{balay2019petsc,
   title={PETSc users manual},
   author={Balay, Satish and Abhyankar, Shrirang and Adams, Mark and Brown, Jed and Brune, Peter and Buschelman, Kris and Dalcin, Lisandro and Dener, Alp and Eijkhout, Victor and Gropp, W and others},
   year={2019},
   institution = {Argonne National Laboratory}
 }

@inproceedings{wkeglarczyk2018kernel,
  title={Kernel density estimation and its application},
  author={W{\k{e}}glarczyk, Stanis{\l}aw},
  booktitle={ITM Web of Conferences},
  volume={23},
  pages={00037},
  year={2018},
  organization={EDP Sciences}
}

@article{Carmona2012ProbabilisticAO,
  title={Probabilistic Analysis of Mean-Field Games},
  author={Ren{\'e} A. Carmona and François Delarue},
  journal={SIAM J. Control. Optim.},
  year={2012},
  volume={51},
  pages={2705-2734}
}

@article{Achdou2012MeanFG,
  title={Mean Field Games: Numerical Methods for the Planning Problem},
  author={Yves Achdou and Fabio Camilli and Italo Capuzzo Dolcetta},
  journal={SIAM J. Control. Optim.},
  year={2012},
  volume={50},
  pages={77-109}
}

@article{Achdou2020MeanFG2,
  title={Mean Field Games of Controls: Finite Difference Approximations},
  journal = {Mathematics in Engineering},
  volume = {3},
  number = {3},
  pages = {1-35},
  year = {2021},
  issn = {2640-3501},
  doi = {10.3934/mine.2021024},
  url = {https://www.aimspress.com/article/doi/10.3934/mine.2021024},
  author = {Yves Achdou and Ziad Kobeissi}
}

@article{nagalur2019first,
  title={First order multi-lane traffic flow model--an incentive based macroscopic model to represent lane change dynamics},
  author={Nagalur Subraveti, Hari Hara Sharan and Knoop, Victor L and van Arem, Bart},
  journal={Transportmetrica B: Transport Dynamics},
  volume={7},
  number={1},
  pages={1758--1779},
  year={2019},
  publisher={Taylor \& Francis}
}

@book{carmona2018probabilistic,
  title={ Probabilistic Theory of Mean Field Games with Applications II Mean Field Games with Common Noise and Master Equations},
  author={Carmona, René and Delarue, Fran{\c{c}}ois},
  year={2018},
  publisher={Springer}
}

@inproceedings{Carmona2018ProbabilisticTO,
  title={Probabilistic Theory of Mean Field Games with Applications I: Mean Field FBSDEs, Control, and Games},
  author={Ren{\'e} A. Carmona and François Delarue},
  year={2018},
  url={https://api.semanticscholar.org/CorpusID:125947789}
}

@article{kamal2014smart,
  title={Smart driving of a vehicle using model predictive control for improving traffic flow},
  author={Kamal, Md Abdus Samad and Imura, Jun-ichi and Hayakawa, Tomohisa and Ohata, Akira and Aihara, Kazuyuki},
  journal={IEEE Transactions on Intelligent Transportation Systems},
  volume={15},
  number={2},
  pages={878--888},
  year={2014},
  publisher={IEEE}
}

@article{Bensoussan1984,
author = {Bensoussan, A. and Frehse, J.},
journal = {Journal für die reine und angewandte Mathematik},
pages = {23-67},
title = {Nonlinear elliptic systems in stochastic game theory.},
url = {http://eudml.org/doc/152630},
volume = {350},
year = {1984},
}

@article{Bensoussan2014MeanFG,
  title={Mean Field Games with a Dominating Player},
  author={Alain Bensoussan and Michael Chau and Phillip S. C. Yam},
  journal={Applied Mathematics \& Optimization},
  year={2014},
  volume={74},
  pages={91-128},
  url={https://api.semanticscholar.org/CorpusID:117798101}
}

@article{McDonald2018PreconditioningAI,
  title={Preconditioning and Iterative Solution of All-at-Once Systems for Evolutionary Partial Differential Equations},
  author={Eleanor McDonald and Jennifer Pestana and Andrew J. Wathen},
  journal={SIAM J. Sci. Comput.},
  year={2018},
  volume={40},
  url={https://api.semanticscholar.org/CorpusID:43934903}
}

\appendix

\section{One-class traffic}\label{appendix_one-class}
\begin{figure}[!ht]
     \centering
     \includegraphics[width=\textwidth]{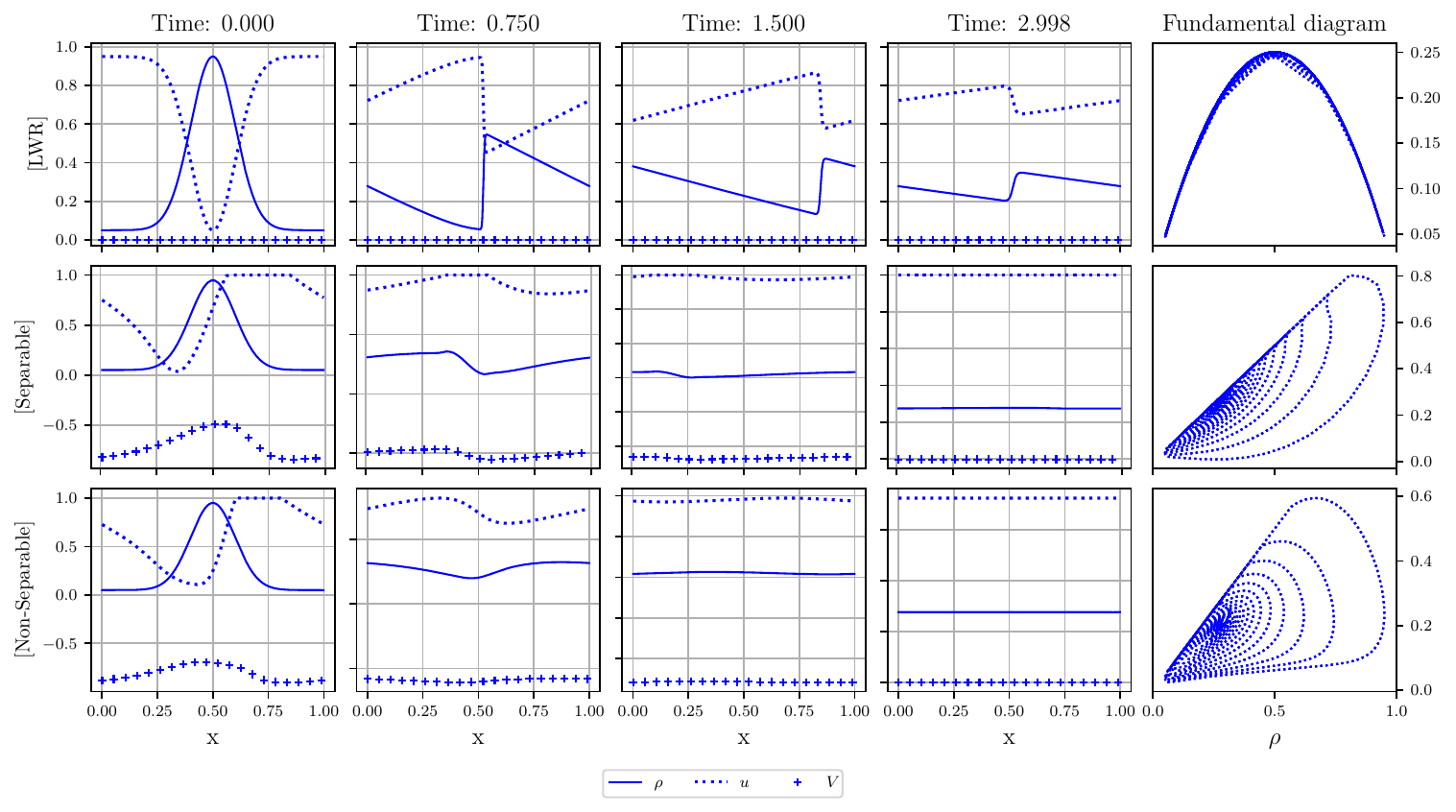}
      \caption{\centering One-class (\textcolor{blue}{cars}): Evolution of density $\rho$, speed $u$, and optimal cost $V$ in four-time values (four columns). The rows show the three cost functions [LWR], [Separable], and [Non-Separable]. The last column shows the fundamental diagram.}
        \label{fig.macro_mfg_1class}
    
\end{figure}

\begin{figure}[!ht]
     \centering
     \includegraphics[width=\textwidth]{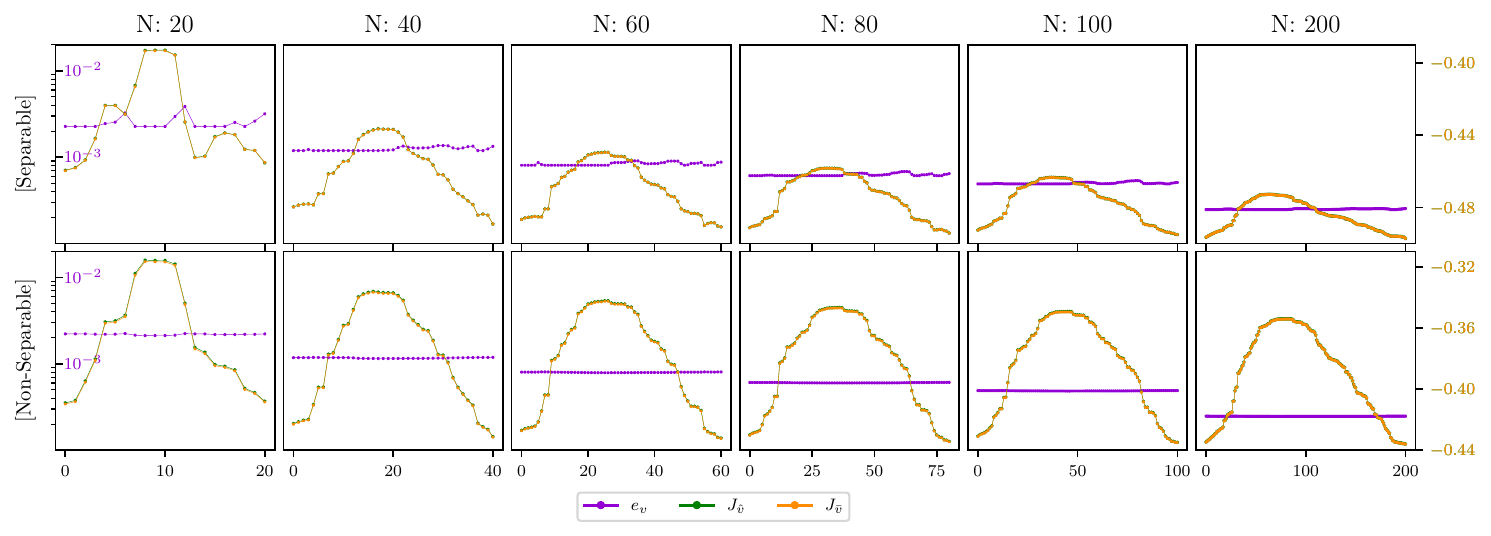}
      \caption{\centering The columns represent the different total number of cars $N$. The rows correspond to the [Separable] and [Non-separable] cost functions. Each sub-figure shows three curves: the $L_{\infty}$ norm \textcolor{darkpurple}{$e_v=\left\|\hat{v}-\Bar{v}\right\|_{\infty}$}, the cost for NMFE-constructed controls \textcolor{darkgreen}{$J_{\hat{v}}=J_{i_{1}}(\hat{v}_{i_{1}},\hat{v}_{-i_{1}})$}, and the cost for the best response strategies \textcolor{darkorange}{$J_{\bar{v}}=J_{i_{1}}(\Bar{v}_{i_{1}},\hat{v}_{-i_{1}})$}, all w.r.t. car's index (idx).
      }
        \label{fig.micro_1class_micro_plots}
    
\end{figure}

For one class of traffic model proposed in \cite{huang_game-theoretic_2020}, we plot in Figure \ref{fig.macro_mfg_1class} the evolution of the density, the velocity field, and the optimal control, as well as the fundamental diagram for the three examples of the cost functions [LWR], [Separable], and [Non-Separable].\\
\\
Figure \ref{fig.micro_1class_micro_plots} shows the $L_{\infty}$ norm $e_v=\left\|\hat{v}-\Bar{v}\right\|_{\infty}$. The cost function when cars move with their NMFE-constructed controls  $J_{\hat{v}}=J_{i_{1}}(\hat{v}_{i_{1}},\hat{v}_{-i_{1}})$, and the cost when cars move with their best response strategies  $J_{\bar{v}}=J_{i_{1}}(\Bar{v}_{i_{1}},\hat{v}_{-i_{1}})$.\\

\begin{figure}[!ht]
     \centering
     \includegraphics[width=0.8\textwidth]{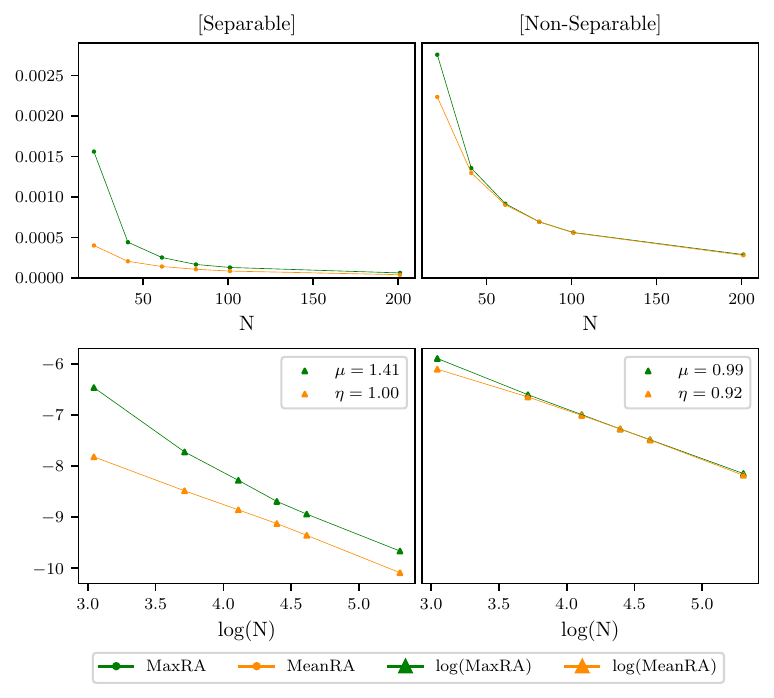}
      \caption{\centering The columns correspond to the [Separable] and [Non-separable] cost functions. The first row shows the mean accuracy (\textcolor{darkorange}{MeanRA}) and maximum accuracy (\textcolor{darkgreen}{MaxRA}) for different total number of cars $N$. The second row shows their respective log-log plots, with $-\mu$ and $-\eta$ being the slopes.}
    \label{fig.micro_1class_RA}
    
\end{figure}

In Figure \ref{fig.micro_1class_RA} we plot in the first row the mean accuracy (MeanRA) and maximum accuracy (MaxRA) with respect to the number of the cars, and in the second row their respective log-log plots.
\newpage

\end{document}